%%%%%%%%%%%%%%%%%%%%%%%%%%%%%%%%%%%%%%%%%%%%%%%%%%%%%%%%%%%%%%%%%%
%%   A Hochschild Homology Euler Characteristic for Circle Actions
%%%%%%%%%%%%%%%%%%%%%%%%%%%%%%%%%%%%%%%%%%%%%%%%%%%%%%%%%%%%%%%%%%
%%   by Ross Geoghegan and Andrew Nicas
%%%%%%%%%%%%%%%%%%%%%%%%%%%%%%%%%%%%%%%%%%%%%%%%%%%%%%%%%%%%%%%%%%
%%   Version:  August 31, 1998 
%%%%%%%%%%%%%%%%%%%%%%%%%%%%%%%%%%%%%%%%%%%%%%%%%%%%%%%%%%%%%%%%%%
%%   TeX file for use with AMSTeX 
%%%%%%%%%%%%%%%%%%%%%%%%%%%%%%%%%%%%%%%%%%%%%%%%%%%%%%%%%%%%%%%%%%

%\input amstex
\documentstyle{amsppt}

\NoRunningHeads
\magnification\magstep1
\pageheight{9truein}
\pagewidth{6.5truein}
\mathsurround=1pt
\NoBlackBoxes

%%%%%%%%%%%%%%%%%%%%%%%%%%
%%%   topmatter.tex       
%%%%%%%%%%%%%%%%%%%%%%%%%%

\topmatter

\title
A Hochschild Homology \hbox{Euler Characteristic} for Circle Actions
\endtitle

\author
Ross Geoghegan$^*$ and Andrew Nicas$^{**}$
\endauthor

\thanks
{
\flushpar
\,\,$^*$Partially supported by the National Science Foundation.
\newline
$^{**}$Partially supported by the Natural Sciences and Engineering
Research Council of Canada.
\newline
The authors acknowledge the hospitality of the Institute for Advanced
Study and support from NSF grant DMS 9304580.
}
\endthanks

\address
Department of Mathematics, SUNY Binghamton,
Binghamton, NY 13902--6000, USA
\endaddress

\email
ross\@math.binghamton.edu
\endemail

\address
Department of Mathematics, McMaster University,
Hamilton, Ontario L8S 4K1, Canada
\endaddress

\email
nicas\@mcmaster.ca 
\endemail

%  Math Subject Classifications

\abstract
We define an ``$S^1$-Euler characteristic'', $\tilde \chi_{S^1}(X)$,
of a circle action on a compact manifold or finite complex $X$.
It lies in the first Hochschild homology group $HH_1(\Bbb Z G)$
where $G$ is the fundamental group of $X$.
This $\tilde \chi_{S^1}(X)$ is analogous in many ways to the ordinary Euler
characteristic.
One application is an intuitively satisfying formula for
the Euler class (integer coefficients) of the normal bundle to a smooth
circle action without fixed points on a manifold.  In the special case of a
$3$-dimensional Seifert fibered space, this formula is particularly
effective.
\endabstract


\endtopmatter

\document

%%%%%%%%%%%%%%%%%%%%%%%%%%
%%%    defs.tex           
%%%%%%%%%%%%%%%%%%%%%%%%%%

%symbols

\def\R{\Bbb {R}}
\def\Z{\Bbb {Z}}
\def\Q{\Bbb {Q}}

\def\F{\Cal {F}}
\def\G{\Cal {G}}
\def\L{\Cal {L}}

\def\vep{\varepsilon}
\def\ra{\rightarrow}
\def\x{\times}
\def\tchi{\tilde \chi}
\def\tChi{\widetilde {\Cal X}}
\def\del{\partial}
\def\tE{\tilde E}
\def\injects{\rightarrowtail}

\def\Hom{\operatorname{Hom}}

\def\euler{\operatorname{Eul}}
\def\trace{\operatorname{trace}}
\def\Fix{\operatorname{Fix}}
\def\id{\operatorname{id}}

\def\DT{\operatorname{DT}}
\def\PD{\operatorname{PD}}
\def\Wh{\operatorname{Wh}}

%chapter labels

\def\chA{1}         %section 1   Algebraic Preliminaries
\def\chN{2}         %section 2   A refined higher Euler characteristic
\def\chF{3}         %section 3   Filtered cell complexes
\def\chS{4}         %section 4   S^1--CW complexes
\def\chV{5}         %section 5   T^2--actions
\def\chT{6}         %section 6   Three dimensional Seifert Fiber Spaces
\def\chG{A}         %appendix A  Some consequences of \tchi'_1(X) \neq 0

%labels for \chA
\def\Alabpa{1}  %Proposition ..g \otimes g^{-1}g_C is a cycle..
\def\Alabpb{2}  %Proposition ..\mu_{\Z}(z) = \sum_i..
\def\Alabda{3}  %Definition ..Dennis trace..
\def\Alabpc{4}  %Proposition ..lies in $HH_1(\Z G)_{C(1)}.

%labels for \chN
\def\Nlabda{1}  %Definition \tChi_1
\def\Nlabpa{2}  %Proposition ..\tChi_1(X,A) becomes a derivation
\def\Nlabpe{3}  %Proposition (change of basis for \tChi_1(X,A))
\def\Nlabdc{4}  %Definition  \tChi'_1(X,A)..\tChi''_1(X,A)
\def\Nlabdb{5}  %Definition  \tchi_1(X,A)..\tchi'_1(X,A)..\tchi''_1(X,A)
\def\Nlabpb{6}  %Proposition (\tchi_1(X,A) = \tchi'_1(X,A) , \tchi''_1(X,A))
\def\Nlabde{7}  %Definition ..\chi_1(X,A), \chi'_1(X,A), \chi''_1(X,A)..
\def\Nlabta{8}  %Theorem Torsion Formula
\def\Nlabla{9}  %Label 
\def\Nlabca{10} %Theorem Simple Homotopy Invariance
\def\Nlabcb{11} %Corollary \tChi''_1(X,A) is invariant under subdivision
\def\Nlabpf{12} %Theorem Homotopy Invariance
\def\Nlabdd{13} %Definition  ..compact polyhedron..
\def\Nlabpd{14} %Proposition For (Z,C) \subset (X,A)..
 %Lemma (combining tensors)

%labels for \chF
\def\Flabda{1}  %Definition ..filtered cell complex..
\def\Flabdd{2}  %Definition ..filtered cell complex X is {\it polyhedral} 
\def\Flabdc{3}  %Definition ..\tChi_1(X; \F)(\gamma)=..
\def\Flabpa{4}  %Proposition T_1(\tilde \del^{n,n-1} .. = \sum..
\def\Flabpb{5}  %Proposition ..\tChi_1(X)(\gamma) = \tChi_1(X; \F)(\gamma)..
\def\Flabta{6}  %Theorem  (main computational tool)

%labels for \chS
\def\Slabda{1}  %Definition ..attaching S^1--cells..
\def\Slabdb{2}  %Definition ..S^1--relative CW complex (X,A)..
\def\Slabpd{3}  %Proposition ..is a filtered cell complex in which..
\def\Slabtd{4}  %Proposition ..S^1--triangulation..
\def\Slabla{5}  %Label \tilde{\partial}^{n,n-1}_{n+1} ...
\def\Slabdc{6}  %Definition ..S^1--Euler characteristic..
\def\Slablb{7}  %Label \tilde D^{\gamma,n,n-1}_n ...
\def\Slabpc{8}  %Proposition \tChi_1(X;\F)(\gamma) = \tchi_{S^1}(X)
\def\Slabta{9}  %Theorem (main theorem)
\def\Slabtb{10} %Theorem ..\chi_1(X)(\gamma) = \chi(X/S^1)\{...
\def\Slabtf{11} %Theorem .. even dimensional .. \chi_1(X)(\gamma)=0.
\def\Slabte{12} %Theorem \PD_X(\Eul(\nu))=..
\def\Slabtg{13} %Theorem  ..\PD(\euler(\nu)) = .. (3 dim Seifert fibered)

%labels for \chV
\def\Vlabta{1}  %Proposition .. T^2--manifold M .. polyhedral filtered ..
\def\Vlabma{2}  %Lemma  \tilde D_*\:C_*(\tilde T^2) .. is given by ..
\def\Vlabtb{3}  %Theorem  \tChi_1(X;\F)(\gamma)=0
\def\Vlabtc{4}  %Theorem  \tChi''_1(X)(\gamma)\circ i = 0 and i^*\tchi'_1(X) =0
\def\Vlabca{5}  %Corollary

%labels for \chT
\def\Tlabpa{1}  %Proposition ..g^i_j is not conjugate to g^k_\ell.
\def\Tlabpb{2}  %Proposition ..C(g^i_j)--.. of \pi''(\tchi_{S^1})..
\def\Tlabtb{3}  %Proposition ..\G(X) is cyclic and..
\def\Tlabca{4}  %Corollary  ..Z(G) is infinite cyclic and..
\def\Tlabpd{5}  %Proposition  (what the invariants detect)
\def\Tlabpc{6}  %Lemma \{\gamma_0\}..has infinite order iff..
\def\Tlabtc{7}  %Theorem  .. features of the Seifert ...
\def\Tlabpe{8}  %Proposition  ... \chi_V(\Sigma) ...

%labels for \chG
\def\Glabta{1}  %Theorem  ..Z(G) is infinite cyclic.
\def\Glabtb{2}  %Theorem  ..\G(X) is cyclic and..
\def\Glabca{3}  %Corollary  ..Z(G) is infinite cyclic and generated..
\def\Glabpa{4}  %Proposition  ..in at most one component.
\def\Glabra{5}  %Remark

%labels for references
\def\AJM{GN$_1$}   %Parametrized Lefschetz-Nielsen Fixed Point..
\def\TTTF{GN$_2$}  %Trace and Torsion in the..
\def\HEC{GN$_3$}   %Higher Euler Characteristics (I)

\baselineskip=20pt plus1pt

%%%%%%%%%%%%%%%%%%%%%%%%%%
%%%   sec00.tex           
%%%%%%%%%%%%%%%%%%%%%%%%%%

\heading
\S 0. Introduction
\endheading

In this paper we introduce
the ``$S^1$--Euler characteristic''
of a finite $S^1$--space,
i.e., of a topological space $X$ with a given circle action,
where $X$ is partitioned as an ``$S^1$--CW complex''
into finitely many ``$S^1$--cells''
(see \S \chS). %xref
We define the
$S^1$--Euler characteristic both geometrically
and as an algebraic trace
and show how it fits into the broad picture of 
algebraic and differential topology and
of algebraic $K$--theory.
A smooth circle action on a compact smooth manifold
is known to admit the structure of a finite
$S^1$--CW complex
\cite{I$_1$, I$_2$} %bibref
and therefore has an $S^1$--Euler characteristic.

Recall that an $S^1$--CW complex $X$ is a union
$X = \bigcup_{n \geq 0} \bigcup_j c^n_j$
of its $S^1$--cells, $c^n_j$;
here the $j$--th $S^1$--$n$--cell
$c^n_j$ is equivariantly modeled on $S^1/H_j \x D^n$
where $H_j \leq S^1$ is the (closed) isotropy subgroup
of the ``interior'' of the $S^1$--cell
(see {\S \chS} %xref
for details).
The isotropy group $H_j$ is either finite cyclic
in which case we define
$|H_j|$ to be the order of $H_j$,
or else $H_j =S^1$ in which case we define $|H_j| =0$.
We will always assume that $X$ is connected.
Write $G$ for the fundamental group of $X$ and let
$\Z G$ denote the integral groupring of $G$.
Our new invariant lies in the Hochschild homology group $HH_1(\Z G)$.
Recall that this abelian group is the first homology of a
certain chain complex;
$1$--chains in this complex are
by definition elements of $\Z G \otimes_{\Z} \Z G$
(see {\S \chA} %xref
for a detailed account).
Thus we write Hochschild $1$--chains as finite sums
$\sum_j n_j u_j \otimes v_j$ where
$u_j, v_j \in G$ and the $n_j$ are integers.

The $S^1$--Euler characteristic,
denoted by $\tchi_{S^1}(X) \in HH_1(\Z G)$,
is defined to be the homology class represented by the
Hochschild $1$--cycle:
$$
\sum_{n\geq 0} (-1)^{n+1}\sum_j \sum^{|H_j|}_{i=1} g_{j,n}
\otimes g^{-1-i}_{j,n}.
\tag{0.1} %xref
$$
Here, the element $g_{j,n} \in G$ (dependent on a choice of basepath)
has, as a representative loop,
an orbit in the interior of the
$S^1$--$n$--cell $c^n_j$ traversed once
(note that during the action of $S^1$, this orbit is
traversed $|H_j|$ times if $|H_j| \neq 0$,
whereas
if $|H_j|=0$ the orbit is trivial and $g_{j,n}=1$).
It turns out that $\tchi_{S^1}(X)$ is independent of basepaths;
see \S \chS. %xref

Formula (0.1) %xref
for $\tchi_{S^1}(X)$ is geometric in the sense that it is described by
certain orbits of the $S^1$--action.
An ``Euler characteristic'' should have an equivalent and more computable
definition as an algebraic trace.
We now give such a description of $\tchi_{S^1}(X)$.

An $S^1$ action on a finite CW complex $X$ defines a homotopy
$X \x I \ra X$ from $\id_X$, the identity map of $X$, to itself.
This homotopy is not necessarily cellular, but it can be deformed
rel $X \x \{0,1\}$ to a cellular homotopy $F\:X \x I \ra X$.
Orient the cells of $X$ and choose  compatibly oriented lifts
in the universal cover, $\tilde X$, of each cell of $X$.
Using cellular chains, 
let $\tilde D\:C_*(\tilde X) \ra C_*(\tilde X)$ be the
chain homotopy corresponding
to the lift $\tilde F\: \tilde X \x I \ra \tilde X$ of
$F$ with $\tilde F_0 = \id_{\tilde X}$.
Let $\tilde \del\:C_*(\tilde X) \ra C_*(\tilde X)$ be the
boundary $\Z G$--homomorphism.
Then we can form the Hochschild $1$--chain (which turns
out to be a $1$--cycle):
$$
\sum_{i,j} \tilde \del_{ij} \otimes \tilde D_{ji}.
\tag{0.2} %xref
$$
Here we regard $\tilde \del = (\tilde \del_{ij})$
and $\tilde D = (\tilde D_{ij})$ as square $\Z G$--matrices;
see Definition \chN.\Nlabda. %xref
One of the main theorems of this paper,
Theorem \chS.\Slabta, %xref
asserts that the homology classes of the 
Hochschild $1$--cycles
(0.1) %xref
and
(0.2) %xref
are essentially the same.
We consider the expression (0.2) %xref
to be a ``trace'' because its homology class generalizes the
``Dennis trace'' of algebraic $K$--theory
(see \S \chN). %xref

One application of the $S^1$--Euler characteristic
is a formula for the Euler class of
the normal bundle to a smooth circle action without fixed
points on a smooth closed oriented manifold $X$.
With the above notation,
the Poincar\'e dual of this Euler class turns out to be
$$
\sum_{n \geq 0} (-1)^n \sum_j \{g_{j,n}\} \in H_1(X; \Z)
\tag{0.3} %xref
$$
(see Theorem \chS.\Slabte). %xref
\footnote{
Notice the striking analogy between 
(0.3) %xref
and the well known identification of the
Poincar\'e dual of the Euler class of the tangent bundle
of a smooth closed oriented manifold $Y$
with the Euler characteristic,
$\sum_{n \geq 0} (-1)^n \sum_j \{y_{j,n}\} \in H_0(Y; \Z) \cong \Z$
(here $y_{j,n}$ is a point, thought of as a $0$--cycle,
in the interior of the $j$--th 
$n$--cell of some finite CW structure on $Y$
so that $\sum_j \{y_{j,n}\}$ is identified
with the number of $n$--cells).  }  %end of footnote 

As a special case, suppose $X$ is an oriented
$3$--dimensional
Seifert fibered space whose
(oriented) quotient surface is $\Sigma$.
The $S^1$--CW structure can then be
chosen so that each singular fiber is an $S^1$--$0$--cell.
Let
$g_{1,0},\dots ,g_{r,0} \in G$ represent the singular fibers and let
$\gamma_0 \in G$ represent an ordinary fiber.
Formula (0.3) %xref
tells us that the
Poincar\'e dual to the Euler class of the normal bundle is
$$
(\chi(\Sigma) - r)\{\gamma_0\} +
\sum_{j=1}^r \{g_{j,0}\}
~\in H_1(X; \Z).
\tag{0.4} %xref
$$
(See Theorem {\chS.\Slabtg}.) %xref
Furthermore, the equivalence of
(0.1) %xref
and
(0.2) %xref
implies that (0.4) %xref
can expressed by the computable formula
$$
A(\sum_{i,j} \tilde \partial_{ij}D_{ji})
\tag{0.5} %xref
$$
where the integer matrix $(D_{ij})$ is the matrix of the chain homotopy
$D\:C_*(X) \ra C_*(X)$ (covered by $\tilde D$) which is induced by the
homotopy $F\:X \x I \ra X$ given by (a cellular approximation to)
the circle action,
and the homomorphism
$A\: \Z G \ra G/[G,G] \equiv G_{\text{ab}} \cong H_1(X;\Z)$
is the unique extension of the abelianization homomorphism
$G \ra G/[G,G]$;
see \cite{GN$_3$}.
\footnote{
The formula (0.4) %xref
is easily obtained for rational coefficients using the fact that a
finite covering of the Seifert fibered space is a circle bundle.
To the best of our knowledge, formulas
(0.4) %xref
and 
(0.5) %xref
with integral coefficients are new.  } %end of footnote

Two further results reinforce the analogy between our
$\tilde \chi_{S^1}(X)$
and the classical Euler characteristic:

\flushpar
$\bullet~$ The classical Euler characteristic of a closed
odd-dimensional manifold is zero.
We show that the image of $\tilde \chi_{S^1}(X)$
under the natural homomorphism $HH_1(\Z G) \ra H_1(G)$ (which
takes a Hochschild homology class represented by 
$\sum_i n_i g_i \otimes h_i$ to $A(\sum_i n_i g_i)$) is zero
whenever $X$ is a closed even-dimensional smooth $S^1$--manifold
(see Theorem {\chS.\Slabtf}  %xref
and its Addendum). %xref

\flushpar
$\bullet~$ The classical Euler characteristic of a finite CW complex
vanishes if the complex admits a circle action with
finite isotropy.
We show in
{\S \chV} %xref
that the $S^1$--Euler characteristic is analogously an
obstruction to the existence of an $S^1 \x S^1$--action
with finite isotropy extending a given $S^1$--action
(Corollary \chV.\Vlabca). %xref

Since $\tchi_{S^1}(X)$
lies in the Hochschild homology group $HH_1(\Z G)$,
the question naturally arises as to whether
it is necessarily contained in the image
of the Dennis trace homomorphism $\DT: K_1(\Z G) \ra HH_1(\Z G)$
(see Definition \chA.\Alabda). %xref
Our extensive calculation of $\tchi_{S^1}(X)$ in case $X$ is
a $3$--dimensional Seifert fibered space (\S \chT) %xref
shows that this is often {\it not\,} true for these examples.

An ``Euler characteristic'' should be topologically invariant.
In the case of the $S^1$--Euler characteristic, the interpretation
of this principle is not immediately clear.
The precise formulation of topological invariance requires the
introduction of some new invariants
(presented here in \S \chN) %xref
of a somewhat technical nature
which extend invariants first developed in \cite{GN$_3$}. %bibref
In fact, $\tchi_{S^1}(X)$ canonically decomposes into two pieces, 
the first of which (after the removal of a slight ambiguity)
is computed from a homotopy invariant of
the underlying space $X$ and the second of which is computed from a
{\it simple} homotopy invariant of $X$.
The connection between
{\S \chN} %xref
and the results discussed in this introduction only becomes clear in
{\S \chS}, %xref
the core section of this paper.
In {\S \chF}, %xref
we prove a key technical result (Theorem \chF.\Flabta) %xref
needed for our discussion of $S^1$--CW complexes
in {\S \chS}. %xref
The $S^1$--Euler characteristic of an $S^1$-complex
is formally introduced in
{\S \chS} %xref
(Definition \chS.\Slabdc). %xref
Its expression in terms of the topological invariants of
{\S \chN} %xref
is given by Theorem \chS.\Slabta,  %xref
and its application to formula (0.3) %xref
for the Euler class of
the normal bundle to a smooth circle action
is given by Theorem {\chS.\Slabte}.  %xref
Formula (0.4) %xref
is the content of Theorem \chS.\Slabtg. %xref
The vanishing theorems for $\tchi_{S^1}(X)$ for a circle action which 
extend to an $S^1 \x S^1$--action with finite isotropy are
proved in {\S \chV}.   %xref 
In {\S \chT} %xref
we compute $\tchi_{S^1}(X)$ when $X$ is an
oriented $3$--dimensional Seifert fibered space
(Theorem \chT.\Tlabpb)
and discuss how much of that structure is detected by $\tchi_{S^1}(X)$
(Theorem \chT.\Tlabtc). %xref
This extended example provides a good illustration of the theory of the
$S^1$--Euler characteristic.
In Appendix {\chG} %xref
we show how the invariants introduced in
{\S \chN} %xref
give higher order analogs of Gottlieb's Theorem \cite{Go}. %bibref

Finally, we remark that the use of tildes in our notation,
as in $\tilde \chi_{S^1}(X)$, is meant to convey to the reader
that the invariant in question is a feature of the universal cover 
$\tilde X$ as opposed to something which can be naively read off from
$X$ itself.

\remark{Acknowledgment}
We are grateful to Matthew Brin for making available his
unpublished lecture notes
\cite{B} %bibref
on Seifert fibered spaces, and
for a number of helpful conversations.
\endremark

%%%%%%%%%%%%%%%%%%%%%%%%%%
%%%   sec01.tex           
%%%%%%%%%%%%%%%%%%%%%%%%%%

\heading
\S \chA. Algebraic Preliminaries
\endheading

In the this section we review
some background material concerning Hochschild
homology, traces and Whitehead torsion;
some general references for these topics are
\cite{C}, %bibref
\cite{I}  %bibref
and
\cite{L}. %bibref

Let $R$ be a commutative ground ring and
let $S$ be an associative $R$--algebra with unit.
If $M$ is an $S$--$S$ bimodule 
(i.e., a left and right $S$--module satisfying 
$(s_1 m) s_2 = s_1 (m s_2)$ for all $m \in M$,
and $s_1, s_2 \in S$),
the {\it Hochschild chain complex} $\{C_* (S,M),d\}$
consists of 
$C_n(S,M) = S^{\otimes n} \otimes M$
where $S^{\otimes n}$  is
the tensor product of $n$ copies of $S$ and
$$
\align
d(s_1\otimes \cdots \otimes s_n \otimes m)&=s_2
\otimes \cdots \otimes s_n \otimes m s_1 \\
& +\sum^{n-1}_{i=1}(-1)^i s_1\otimes \cdots \otimes
s_i s_{i+1}\otimes \cdots \otimes s_n\otimes m \\
&+(-1)^n s_1\otimes \cdots \otimes s_{n-1}\otimes s_n m.
\endalign
$$
The tensor products are taken over $R$. 
The $n$--th homology of this complex is the
$n$--th
{\it Hochschild homology of} $S$ {\it with coefficient bimodule} $M$.  
It is denoted by $HH_n(S,M)$.
If $M=S$ with the standard $S$--$S$ bimodule structure
then we usually write $HH_n(S)$ for $HH_n(S,M)$.

We will be concerned mainly with $HH_1$ and $HH_0$ which
are computed from 
$$
\alignat2
\cdots \longrightarrow &S\otimes S \otimes
M \ \ \ \ {\overset d \to \longrightarrow}
&&\qquad S \otimes M  \ \ \ \ {\overset d \to \longrightarrow} \qquad M \\ 
&s_1\otimes s_2\otimes m \ \ \ \ \mapsto &&\qquad s_2 \otimes m
s_1 - s_1 s_2 \otimes m + s_1 \otimes s_2 m \\
& &&\qquad s \otimes m \ \ \ \  \mapsto \qquad m s -s m
\endalignat
$$

Next, we consider traces in Hochschild homology.
If $A$ is a square matrix over $M$, we interpret its trace
$\sum_i A_{ii}$ as an element of $M$
(i.e., as a Hochschild $0$--cycle).
The corresponding homology class is denoted by
$T_0(A) \in HH_0(S,M)$.
If $A^i$, $i=1, \ldots, n$,
are $q_i \times q_{i+1}$ matrices
over $S$ and $B$ is a $q_{n+1} \times q_1$ matrix over $M$,
we define
$A^1 \otimes \cdots \otimes A^n \otimes B$
to be the $q_1 \x q_1$ matrix with entries
in the $R$--module $S^{\otimes n} \otimes M$ given by 
$$
(A^1 \otimes \cdots \otimes A^n \otimes B)_{ij} =
\sum_{k_2, \ldots, k_{n+1}}
A^1_{i, k_2} \otimes A^2_{k_2, k_3} \otimes \cdots \otimes
A^n_{k_n, k_{n+1}} \otimes B_{k_{n+1}, j}.
$$
The {\it trace} of
$A^1 \otimes \cdots \otimes A^n \otimes B$,
written
$\trace(A^1 \otimes \cdots \otimes A^n \otimes B)$,
is
$$
\sum_{k_1,k_2, \ldots, k_{n+1}}
A^1_{k_1, k_2} \otimes A^2_{k_2, k_3} \otimes \cdots \otimes
A^n_{k_n, k_{n+1}} \otimes B_{k_{n+1}, k_1}.
$$
which we interpret as a Hochschild $n$--chain.
Observe that the $1$--chain
$\trace(A \otimes B)$ is a cycle if and only if
$\trace(AB)=\trace(BA)$,
in which case we denote its homology class by
$T_1(A \otimes B)\in HH_1(S,M)$.
In this paper $S$ will usually be a groupring
over the ground ring $R$ and $M=S$.

We will use the notation $G_1$ for the set of
conjugacy classes of a group $G$, and $C(g)$
for the conjugacy class of $g \in G$.
The partition of $G$ into the union of its conjugacy classes
induces a  direct sum decomposition of
$HH_*(\Z G) \equiv HH_*(\Z G, \Z G)$ as follows:
each generating chain
$c=g_1 \otimes \cdots \otimes g_n \otimes m$
can be written in {\it canonical form} as
\linebreak
$ g_1 \otimes \cdots \otimes g_n
\otimes g_n^{-1} \cdots g_1^{-1}  g$
where we think of $g = g_1 \cdots g_n m \in G$ as 
``marking'' the conjugacy class $C(g)$.
All the generating chains occurring in the boundary $d(c)$
are easily seen to
have markers in $C(g)$ when put into canonical form.
For $C \in G_1$ let $C_*(\Z G)_C$ be the subgroup of
$C_*(\Z G) \equiv C_*(\Z G,\Z G)$
generated by those generating  chains whose
markers lie in $C$.
The decomposition
$$
\Z G  \cong \bigoplus_{C \in G_{1} } \Z C
$$
as a direct sum of abelian groups determines a decomposition
of chain complexes
$$
C_*(\Z G) \cong 
\bigoplus_{C \in G_{1} } C_*(\Z G)_C.
$$
There results a natural isomorphism
$HH_*(\Z G) \cong 
\bigoplus_{C \in G_{1} } HH_*(\Z G)_C$ where
the summand $HH_*(\Z G)_C$ corresponds to the homology
classes of Hochschild cycles marked by the elements of $C$.
We call this summand the {\it $C$--component}.

Given any $\Z G$--$\Z G$  bimodule $N$ let
$\overline N$ be the left $\Z G$ module whose underlying abelian
group is $N$ and whose left module structure is given by
$g m = g \cdot m \cdot g^{-1}$.
There is a natural isomorphism
$\mu_N\:HH_*(\Z G,N) @>{\cong}>> H_*(G, \overline N)$
which is induced from an isomorphism of the Hochschild
complex to the bar complex for computing group homology;
see \cite{I, Theorem 1.d}.  %bibref
The decomposition
$\overline {\Z G} \cong \bigoplus_{C \in G_{1} } \Z C$
is a direct sum of left $\Z G$ modules,
inducing a direct sum decomposition
$H_*(G, \overline {\Z G}) \cong \bigoplus_{C \in G_{1} } H_*(G,\Z C)$.
Choosing representatives  $g_C \in C$ we have an
isomorphism of left
$\Z G$ modules   $\Z C \cong  \Z (G/Z(g_C))$ where 
$Z(h)=\{ g \in G ~|~ h = g  h  g^{-1}  \}$  denotes the 
centralizer of $h \in G$.
Since $H_*(G, \Z (G/Z(g_C)))$ is naturally isomorphic
to $H_*(Z(g_C))$,
we obtain a natural isomorphism
$HH_*(\Z G )  \cong \bigoplus_{C \in G_{1} } H_*(Z(g_C))$;
$HH_*(\Z G )_C$ corresponds to the
summand $H_*(Z(g_C))$ under this identification.
In particular $HH_0(\Z G)  \cong \Z G_{1}$,
the free abelian group generated by the conjugacy classes,
and $HH_1(\Z G )  \cong  \bigoplus_{C \in G_{1} } H_1(Z(g_C))$,
the direct sum of the abelianizations of the centralizers.

\proclaim{Proposition \chA.\Alabpa} %xref
The chain $g \otimes g^{-1}g_C$ is a cycle if and only if
$g \in Z(g_C)$.
For $g \in Z(g_C)$,
the homology class of $g \otimes g^{-1}g_C$ 
in $HH_1(\Z G)$ corresponds
to $\{g\} \in H_1(Z(g_C))$. \hfill{ }\qed
\endproclaim

The augmentation $\vep:\Z G \ra \Z$ can be viewed as a
morphism of $\Z G$--$\Z G$ bimodules,
where $\Z$ is given the trivial bimodule structure,
or as a morphism
$\vep:\overline {\Z G} \ra \overline  \Z$ of
left $\Z G$--modules.
Then there is an induced chain map
$C_*(\Z G,\Z G) @>{\vep}>> C_*(\Z G,\Z)$
and a commutative diagram:
$$
\CD
HH_*(\Z G, \Z G )  @>{\vep}>>  HH_*(\Z G, \Z) \\
@V{\mu_{\Z G}}VV   @V{\mu_{\Z}}VV \\
H_*(G, \overline {\Z G} )  @>{\vep}>>  H_*(G, \overline \Z)
\endCD
$$
where the vertical arrows are isomorphisms.
Let $\vep_*\:HH_*(\Z G) \ra H_*(G)$ denote the resulting
homomorphism.

The ``abelianization'' homomorphism
$A\:G \ra G_{\text{ab}} = H_1(G)$
extends to a homomorphism of abelian groups
$A:\Z G \ra G_{\text{ab}}$.
This occurs in the following computation of $\mu_{\Z}$
(\cite{\HEC, Proposition 2.1}): %bibref

\proclaim{Proposition \chA.\Alabpb} %xref
If $\sum_i c_i \otimes n_i \in C_1(\Z G, \Z)$ is
a Hochschild $1$--cycle representing
\linebreak
$z \in HH_1(\Z G, \Z)$,
where $c_i \in \Z G$ and $n_i \in \Z$,
then $\mu_{\Z}(z) = \sum_i A(c_i n_i) \in H_1(G)$.\hfill{ }\qed
\endproclaim

We recall the definition of $K_1$ of a ring.
Let $GL(n,R)$ denote the general linear group
consisting of all $n \x n$ invertible matrices
over $R$, and let $GL(R)$ be the direct limit
of the sequence $GL(1,R) \subset GL(2,R) \subset \cdots$.
A matrix in $GL(R)$ is called {\it elementary} if
coincides with the identity except for a single
off-diagonal entry.
The subgroup $E(R) \subset GL(R)$ generated by
the elementary matrices is precisely the commutator
subgroup of $GL(R)$, and the abelian
quotient group $GL(R)/E(R)$ is, by definition, $K_1(R)$.

\definition{Definition \chA.\Alabda} %xref
The {\it Dennis trace} homomorphism
$
\DT\:K_1(R) @>{}>> HH_1(R)
$
is defined as follows.
If $\alpha \in K_1(R)$ is represented by
an invertible $n \x n$ matrix $A$ then\linebreak
$\DT(\alpha) = T_1(A \otimes A^{-1})~$
(see \cite{I, Chapter 1}).  %bibref
\enddefinition

In case $R=\Z G$, let $\pm G \subset GL(1,\Z G)$ be the
subgroup consisting of $1 \x 1$ matrices of the form
$[\pm g]$, $g \in G$.
The cokernel of the natural homomorphism
$\pm G \ra K_1(\Z G)$
is called the
{\it Whitehead group of $G$}
and is denoted by $\Wh_1(G)$.

\proclaim{Proposition \chA.\Alabpc} %xref
The image of the composite homomorphism:
$$
\pm G @>{i}>> K_1(\Z G) @>{\DT}>> HH_1(\Z G)
$$
lies in $HH_1(\Z G)_{C(1)}$.
\endproclaim

\demo{Proof}
For $g \in G$,
$\DT(i(\pm g)) = \text{homology class of }
g \otimes g^{-1} \in HH_1(\Z G)_{C(1)}$.\hfill{ }\qed
\enddemo

We briefly recall how the torsion of an
acyclic complex
is defined (see \cite{C}). %bibref
Declare two bases of a finitely generated free
right $\Z G$--module to
be {\it equivalent} if the change of basis matrix represents
an element of $K_1(\Z G)$ which lies in the image
of $\pm G @>{i}>> K_1(\Z G)$.
Let $(C, \partial)$ 
be a finitely generated chain complex of right
modules over $\Z G$
such that each $C_i$ is free with a given 
equivalence class of bases.
Suppose that $C$ is acyclic.
Let $\delta\: C \ra C$ be a chain contraction,
$C_{\text{odd}} \equiv \oplus_{i \text{ odd}} C_i$ and
$C_{\text{even}} \equiv \oplus_{i \text{ even}} C_i$.
The restriction of $\partial + \delta$ to
$C_{\text{odd}}$ is an isomorphism
$C_{\text{odd}} \ra C_{\text{even}}$ and so its
matrix with respect to bases chosen
from each of the given equivalence classes
defines an element of $K_1(\Z G)$.
The image of this element in $\Wh_1(G)$,
denoted $\tau(C)$,
is independent of the choice of representatives
of the equivalence classes of bases;
it is called the {\it Whitehead torsion} of $(C, \partial)$.

%%%%%%%%%%%%%%%%%%%%%%%%%%
%%%   sec02.tex           
%%%%%%%%%%%%%%%%%%%%%%%%%%

\heading
\S \chN. A refined higher Euler characteristic
\endheading

The goal of this section is to develop topologically invariant
refinements of the higher Euler characteristics first introduced
in \cite{\HEC}; %bibref
prior knowledge of \cite{\HEC} %bibref
is not needed.
In \S \chS, %xref
these refined invariants will used be prove our main results about
the $S^1$--Euler characteristic.
In \S \chN(A)  %xref
the invariants are defined and their basic properties are established.
Homotopy and simple homotopy invariance is discussed in \S \chN(B). %xref
Their geometric interpretation is given in \S \chN(C) %xref
(which might usefully be read in conjunction with \S \chN(A) %xref
to provide motivation).

\subheading{(A) Definition of the invariants}

Let $(X,A)$ be a topological pair where $X$ is a
$k$--space (i.e., a Hausdorff space which has the
weak topology determined by the family of its 
compact subspaces) and $A \subset X$ is closed.
Let $X^X$ denote the function space of continuous maps
$X \ra X$ with the compact-open topology and let
$(X,A)^{(X,A)} \subset X^X$ be the subspace consisting of
maps of pairs $(X,A) \ra (X,A)$.
Define $\Gamma_{(X,A)} \equiv \pi_1( (X,A)^{(X,A)}, \id)$.
Note that $\Gamma_{(X,A)}$ is abelian because
$(X,A)^{(X,A)}$ is an H-space.
An element
$\gamma \in \Gamma_{(X,A)}$ is represented by a continuous
map of pairs
$F^\gamma \:(X,A) \x I \ra (X,A)$ such that
$F_0^\gamma = F_1^\gamma = \id$.
When $A=\emptyset$ we write $\Gamma_X$ for $\Gamma_{(X, \emptyset)}$.
The inclusion $(X,A)^{(X,A)} \subset X^X$ induces a homomorphism
$\Gamma_{(X,A)} \ra \Gamma_X$.

Suppose $X$ is path connected.
Let $x_0 \in X$ be a basepoint and let $G \equiv \pi_1(X,x_0)$.
Evaluation at $x_0$ yields a map $\eta\:X^X \ra X$.
The image of $\eta_{\#}\:\Gamma_X \ra G$, denoted $\G (X)$,
is called the
{\it Gottlieb subgroup of $G$}.
Let $Z(G)$ denote the center of $G$.
It was shown in
\cite{Got} %bibref 
that $\G (X) \leq Z(G)$ and that if $X$ is also aspherical
then $\G(X)=Z(G)$ and
$\eta_{\#}\:\Gamma_X \ra Z(G)$ is an isomorphism.
In particular,
we may regard $\G (X) \subset G_1$  (where $G_1$ is the set
of conjugacy classes of elements of $G$).
Define
$$
\align
HH_*(\Z G)' &\equiv \bigoplus_{C \in \G (X)} HH_*(\Z G)_C \\
HH_*(\Z G)'' &\equiv \bigoplus_{C \in G_{1} - \G (X) } HH_*(\Z G)_C
\endalign
$$

There is a left action of $Z(G)$ on $HH_*(\Z G)$.
At the level of chains it is defined by
$$
\omega \cdot \big( g_1 \otimes \cdots \otimes g_n \otimes m \big)
~=~
g_1 \otimes \cdots \otimes g_n \otimes (m\omega^{-1})
$$
where $\omega\in Z(G)$.
One easily checks that this action is compatible with the
Hochschild boundary
$d$ and hence makes $HH_*(\Z G)$ into a left
$Z(G)$--module.
The summand $HH_*(\Z G)_C$ is taken by the left action
of $\omega$ isomorphically onto the summand
$HH_* (\Z G)_{C \omega^{-1}}$ where
$C \omega^{-1}$ 
is the conjugacy class $\{g\omega^{-1}~|~g\in C\}$.
Since $\eta_{\#}$ maps $\Gamma_{X}$ into $Z(G)$,
$\eta_{\#}$ defines a left action of
$\Gamma_{X}$
(and thus also a left action of $\Gamma_{(X,A)}$
via the natural homomorphism $\Gamma_{(X,A)} \ra \Gamma_X$)
on $C_*(\Z G, \Z G)$
and on $HH_1(\Z G)$.
The group $HH_*(\Z G)$ splits as a direct sum of
$\Gamma_{(X,A)}$--modules:
$$
HH_*(\Z G) = HH_*(\Z G)' \oplus HH_*(\Z G)''.
$$

Let
$$
\pi'\:HH_*(\Z G) \ra HH_*(\Z G)'
\qquad \text{and} \qquad
 \pi''\:HH_*(\Z G) \ra HH_*(\Z G)''
$$
denote the projections. 
Since the centralizer of any element in  $\G(X)$ is all of $G$,
there is a natural isomorphism of $\Gamma_{(X,A)}$--modules
$HH_*(\Z G)' \cong H_1(G) \otimes_{\Z} \Z \G(X)$.

Given a relative CW complex $(X,A)$, we denote its
relative $n$--skeleton by $(X,A)^n$.
Recall
(see \cite{Sp}) %bibref
that the cellular chain complex $C_*(X,A)$
of $(X,A)$
is defined by
$$
C_j(X,A) = H_j((X,A)^j, (X,A)^{j-1})
$$
with the
boundary homomorphism arising from the long
exact sequence associated to the triple
$((X,A)^j, (X,A)^{j-1}, (X,A)^{j-2})$.
Orient the cells of $(X,A)$,
thus establishing a preferred basis for the chain complex
$(C_*(X,A),\partial)$.

For the rest of this section we assume
that $(X,A)$ is a finite relative CW complex such that:
\roster
\item
$A$ is a $k$--space.

\item
$X$ is path connected.

\item
$X$ has a simply connected universal covering space, $\tilde X$.
\endroster
Define $\tilde A \equiv p^{-1}(A)$ where is $p\:\tilde X \ra X$
is the universal covering projection.
The pair $(\tilde X, \tilde A)$ is a relative CW complex
with $j$--skeleton
$(\tilde X, \tilde A)^j = p^{-1}((X, A)^j)$.
Choose a lift, $\tilde e$,
to the universal cover, $\tilde X$, for each cell $e$ of $(X,A)$,
and orient $\tilde e$ compatibly with $e$.
Regard the cellular chain complex
$(C_*(\tilde X, \tilde A),\tilde \partial)$ as a free right
$\Z G$--module chain complex with preferred basis
$\{ \tilde e \}$;
see \cite{\AJM, \S 1(B)}. %bibref

By the Cellular Approximation Theorem,
each $\gamma \in \Gamma_{(X,A)}$ can be represented by a cellular
homotopy $F^\gamma\: (X,A) \x I \ra (X, A)$ such that
$F^\gamma_0 = F^\gamma_1 = \id _X$.
There is a unique lift,
$\tilde F^\gamma\:(\tilde X, \tilde A) \x I \ra (\tilde X, \tilde A)$,
of $F^\gamma$ such that $\tilde F^\gamma_0 = \id_{\tilde X}$.
Note that $\tilde F^\gamma_1$ is the covering transformation
corresponding to the element $\eta_{\#}(\gamma) \in G$.
Let
$\tilde D^\gamma_*:C_*(\tilde X, \tilde A) \ra C_{*+1}(\tilde X, \tilde A)$
be the chain homotopy induced by $\tilde F^\gamma$.

\example{Sign Convention}
If $\tilde e$ is an oriented $k$--cell of $(\tilde X, \tilde A)$ then
$\tilde D^\gamma_k(\tilde e)$ is the  $(k+ 1)$--chain
\linebreak
$(-1)^{k+1}\tilde F^\gamma_*(\tilde e \times I)\in C_{k+1}(\tilde X, \tilde A)$,
where $\tilde e \times I$ is given the product orientation.
This is consistent with the convention that if $E_{i,\epsilon}$ is the face
of the cube $I^n = [0,1]^n$ obtained by holding the $i^{\text{th}}$
coordinate fixed at $\epsilon = 0$ or $1$,
then the incidence number $[I^n : E_{i,\epsilon}]$ is $(-1)^{i+\epsilon}$.
At the level of cellular
$n$--chains, we have
$\partial_nI^n =  \sum_{i,\epsilon} [I^n : E_{i,\epsilon}]E_{i,\epsilon}$.    
\endexample

Write
$\tilde \partial = \bigoplus_k \tilde \partial_k$,
$\tilde D^\gamma = \bigoplus_k (-1)^{k+1} \tilde D^\gamma_k$
and
$\tilde I = \bigoplus_k (-1)^k \id_k$
(viewed as matrices).
The chain homotopy relation becomes
$\tilde D^\gamma \tilde \partial - \tilde \partial \tilde D^\gamma
 = \tilde I (1- \eta_{\#}(\gamma)^{-1})$.
[Explanation: the minus sign occurs
on the left because of the sign convention
built into the matrix $\tilde D^\gamma$;
the right hand side is thus because the
$0$--end of the homotopy $F^\gamma$ is lifted to the identity,
while the $1$--end is lifted to the covering translation
corresponding to $\eta_{\#}(\gamma)$;
the inversion occurs because we have $G$
acting on the right.]

The proofs of
Propositions 2.5 and 2.2 of \cite{\HEC} %bibref
carry over directly to show, respectively, that
the Hochschild $1$--chain
$\trace(\tilde \partial \otimes \tilde D^\gamma)$ is a cycle
and that the homology class of this cycle,
$T_1(\tilde \partial \otimes \tilde D^\gamma)$,
does not depend on the
choice of the cellular map $F^\gamma$ representing
$\gamma \in \Gamma_{(X,A)}$.

\definition{Definition \chN.\Nlabda} %xref
Define the function $\tChi_1(X,A)\:\Gamma_{(X,A)} \ra HH_1(\Z G)$
by:
$$
\tChi_1(X,A)(\gamma) ~\equiv~ T_1(\tilde \partial \otimes \tilde D^\gamma).
$$
\enddefinition
In case $A = \emptyset$, $\tChi_1(X,\emptyset) \equiv \tChi_1(X)$
is identical to the ``lift of $\chi_1(X)$'' defined in
\linebreak
\cite{\HEC, \S2}. %bibref
By considering lifts of homotopies, we
easily obtain:

\proclaim{Proposition \chN.\Nlabpa} %xref
When $HH_1(\Z G)$ is regarded as a left
$\Gamma_{(X,A)}$--module,
$\tChi_1(X,A)$ becomes a derivation;
i.e.,
$\tChi_1(X,A)(\gamma_1\gamma_2) =
\tChi_1(X,A)(\gamma_1) +
\gamma_1 \cdot \tilde {\Cal X}_1(X,A)(\gamma_2)$.\hfill{ }\qed
\endproclaim

The dependence of $\tChi_1(X,A)$ on the choice of oriented
cell lifts can be made explicit as follows.
Suppose $\L_1 =\{\tilde e^n_j \}$ is a choice of oriented cell lifts
(which we regard as an ordered basis for $C_*(\tilde X, \tilde A)$).
Then any other choice, $\L_2$, of oriented cell lifts is of the
form  $\L_2 =\{\tilde e^n_j (\mu_{n,j} g_{n,j}) \}$
where $g_{n,j} \in G$ and $\mu_{n,j} = \pm 1$.
Let $\tChi_1(X,A;{\L_i})$ be $\tChi_1(X,A)$ computed with respect
to $\L_i$, $i=1,2$.
Let $U$ be the diagonal matrix whose diagonal entries
are $\mu_{n,j} g_{n,j} \in \pm G \subset \Z G$
(this is the ``change of basis matrix'').

\proclaim{Proposition \chN.\Nlabpe} %xref
For $\gamma \in \Gamma_{(X,A)}$,
$$
\tChi_1(X,A;{\L_1})(\gamma) - \tChi_1(X,A;{\L_2})(\gamma) =
(1 - \gamma) \cdot T_1(U \otimes U^{-1}).
$$
Furthermore,
$(1 - \gamma) \cdot T_1(U \otimes U^{-1}) \in HH_1(\Z G)'$.
\endproclaim

\demo{Proof}
By the ``change of basis formula'',
\cite{\AJM, Proposition 3.3}: %bibref
$$
\tChi_1(X,A; {\L_1})(\gamma) - \tChi_1(X,A; {\L_2})(\gamma) =
T_1(U \otimes U^{-1}(1 - \eta_{\#}(\gamma)^{-1})).
$$
Since $\eta_{\#}(\gamma) \in \G (X) \subset Z(G)$,
$$
T_1(U \otimes U^{-1}(1 - \eta_{\#}(\gamma)^{-1})) =
(1 - \gamma) \cdot T_1(U \otimes U^{-1}).
$$
We have
$
\trace(U \otimes U^{-1}) = \sum_{n,j}
g_{n,j} \otimes g^{-1}_{n,j} \in C_1(\Z G, \Z G)_{C(1)}
$
and so
$
(1 - \gamma) \cdot T_1(U \otimes U^{-1})  \in HH_1(\Z G)'
$.\hfill{ }\qed
\enddemo

\definition{Definition \chN.\Nlabdc} %xref
We define derivations:
$$
\tChi'_1(X,A)\:\Gamma_{(X,A)} \ra HH_1(\Z G)'
\qquad \text{and} \qquad
\tChi''_1(X,A)\:\Gamma_{(X,A)} \ra HH_1(\Z G)''
$$
by
$\tChi'_1(X,A) \equiv \pi' \circ  \tChi_1(X,A)$ and
$\tChi''_1(X,A) \equiv \pi'' \circ \tChi_1(X,A)$.
\enddefinition
The fact that $\tChi'_1(X,A)$ and $\tChi''_1(X,A)$ are 
derivations follows from
Proposition {\chN.\Nlabpa}.  %xref
Of course $\tChi_1(X,A) =  (\tChi'_1(X,A), \tChi''_1(X,A))$.

\remark{Remark}
Proposition {\chN.\Nlabpe} %xref
implies that  $\tChi''_1(X,A)$ does not depend on the
choice of oriented cell lifts because
$\pi''((1 - \gamma) \cdot T_1(U \otimes U^{-1})) =0$.
A stronger conclusion is found in
Theorem {\chN.\Nlabca} %xref
below.
\endremark

If $M$ is a left $\Gamma_{(X,A)}$--module then the
cohomology group $H^1(\Gamma_{(X,A)},M)$ is naturally
identified with the quotient,
$\operatorname{Der}(\Gamma_{(X,A)},M) / 
 \operatorname{Inn}(\Gamma_{(X,A)},M)$,
of derivations modulo inner derivations
(recall that a derivation $\Gamma_{(X,A)} \ra M$ is
inner if it is of the form $\gamma \mapsto (1-\gamma)m$
for some $m \in M$).
We denote the image of a derivation
$\Theta\:\Gamma_{(X,A)} \ra M$
in $H^1(\Gamma_{(X,A)},M)$ by $[\Theta]$.

\definition{Definition \chN.\Nlabdb} %xref
We define cohomology classes:
$$
\align
\tchi_1(X,A)  &\equiv [\tChi_1(X,A)] \in  H^1(\Gamma_{(X,A)}, HH_1(\Z G)) \\
\tchi'_1(X,A) &\equiv [\tChi'_1(X,A)] \in  H^1(\Gamma_{(X,A)}, HH_1(\Z G)') 
 \cong  H^1(\Gamma_{(X,A)}, H_1(G) \otimes_{\Z} \Z \G(X))\\
\tchi''_1(X,A) &\equiv [\tChi''_1(X,A)] \in  H^1(\Gamma_{(X,A)}, HH_1(\Z G)'').
\endalign
$$
\enddefinition

When $A = \emptyset$, $\tchi_1(X,\emptyset) \equiv \tchi_1(X)$
is the invariant $\tchi_1(X)$ of
\cite{\HEC, Proposition 2.7}. %bibref
   From the definitions, it is clear that:

\proclaim{Proposition \chN.\Nlabpb} %xref
$\tchi_1(X,A) = (\tchi'_1(X,A), \tchi''_1(X,A))$.\hfill{ }\qed
\endproclaim

\remark{Remark}
In the non-aspherical case our invariants can be non-zero even
when $\chi(X) \neq 0$ (implying $\G(X)$ is trivial).
An example is the real projective plane $P^2$.
It is shown in \cite{\HEC, \S 3(D)} %xref
that $\tchi'_1(P^2) \neq 0$ and $\tchi''_1(P^2) \neq 0$.
\endremark

Some group theoretic and topological consequences of the hypothesis
$\tchi'_1(X) \neq 0$ are given in
Appendix \chG. %xref

\definition{Definition \chN.\Nlabde} %xref
The natural homomorphism $\vep_*\:HH_1(\Z G) \ra H_1(G)$
(see \S 1) %xref
induces a homomorphism:
$$
\vep_*\:H^1(\Gamma_{(X,A)}, HH_1(\Z G)) @>{}>>
H^1(\Gamma_{(X,A)}, H_1(G)) \cong \Hom(\Gamma_{(X,A)}, H_1(G)).
$$
Define three homomorphisms $\Gamma_{(X,A)} \ra  H_1(G) = G_{\text{ab}}$:
$$
\chi_1(X,A)    \equiv \vep_*(\tchi_1(X,A)), \quad
\chi'_1(X,A)   \equiv \vep_*(\tchi'_1(X,A)), \quad
\chi''_1(X,A)  \equiv \vep_*(\tchi''_1(X,A)).
$$
Note that $\chi_1(X,A)  = \chi'_1(X,A) + \chi''_1(X,A)$ by
Proposition {\chN.\Nlabpb}. %xref
\enddefinition
   From
Definition {\chN.\Nlabda} %xref
and
Proposition {\chA.\Alabpb}, %xref
we see that
$$
\chi_1(X,A)(\gamma) = \sum_{k \geq 0} (-1)^{k+1}
A(\trace ( \tilde \partial_{k+1}~ D^\gamma_k))
$$
where the integer matrix $D^\gamma_k$ is the matrix of the chain homotopy
$C_k(X,A) \ra C_{k+1}(X,A)$ induced by a cellular homotopy
$F^\gamma\:(X,A) \x I \ra (X,A)$ representing $\gamma \in \Gamma_{(X,A)}$.
In case $A = \emptyset$, this is 
Definition A$_1$ of \cite{\HEC}. %bibref
In \cite{\HEC} %bibref
we also give very different formulations of $\chi_1(X)$ in terms
of cap products and of stable homotopy.

\subheading{(B) Homotopy and simple homotopy invariance}

Suppose that $(Y,B)$ is a finite relative CW complex such that
$B$ is a $k$--space,
$Y$ is path connected and has a simply connected universal
covering space $\tilde Y$.
Let
\hbox{$h\:(X,A) \ra (Y,B)$}  %WARNING  TeX trickery
be a homotopy equivalence of pairs.
Choose a
homotopy inverse $g\:(Y,B) \ra (X,A)$ for $h$ and
let $H\:(Y,B) \x I \ra (Y,B)$ be a homotopy
$hg \simeq \id_Y$
giving a basepath
\linebreak
$\hat H\:I \ra (Y,B)^{(Y,B)}$ from $hg$ to $\id_Y$.
The map
$$
(X,A)^{(X,A)} \ra (Y,B)^{(Y,B)} \qquad
f \mapsto h \circ f \circ g
$$
together with the basepath $\hat H$ induce an isomorphism 
$h_*\:\Gamma_{(X,A)} \ra \Gamma_{(Y,B)}$.
Since $\Gamma_{(Y,B)}$ is abelian,
$h_*$ is independent of choice of homotopy inverse $g$ and
basepath $\hat H$.

Choose $y_0 = h(x_0) \in B$ as the basepoint for $Y$ and
let $K = \pi_1(Y,y_0)$.
Then $h$ induces an isomorphism $h_{\#}\:G \ra K$.
Recall that the {\it torsion of $h$},  $\tau(h) \in \Wh_1(K)$,
is defined as follows.
Let $h'\:(X,A) \ra (Y,B)$ be a cellular map which is homotopic
to $h$ relative to $A$.
Form the mapping cylinder of $h'$
$$
M(h') \equiv   X \x I \cup_{h'} Y \equiv  (X \x I \sqcup Y)/\thicksim
$$
where ``$\sqcup$'' is disjoint union and $\thicksim$ is the
equivalence relation generated by $(x,1) \thicksim h'(x)$ for $x \in X$.
Let $E \equiv A \x I \cup_{h'|_A} B \subset M(h')$.
The relative CW complex structures on $(X,A)$ and $(Y,B)$ determine
a relative CW complex structure on $(M(h'), E)$;
furthermore, the map $i\:X \ra M(h')$ given by $x \mapsto (x,0)$
is the inclusion of a subcomplex.
Let $K'= \pi_1(M(h'), {\ssize (x_0,0)} )$.
The ``collapse'' map $c\:M(h') \ra Y$ is homotopy equivalence;
in particular, it induces an isomorphism $c_{\#}\:K' \ra K$.
Let $\widetilde{M}(h')$ be the universal cover of $M(h')$ and
let $\tilde E$ be the inverse image of $E$ in $\widetilde{M}(h')$.
Choose oriented lifts of the cells of $M(h')$
to $\widetilde{M}(h')$ (compatible with choice of oriented lifts of
cells of $X$).
Since
$i_*\:C_*(\tilde X, \tilde A) \injects C_*(\widetilde{M}(h'), \tilde E)$,
is a chain homotopy equivalence,
the right $\Z K'$--module complex
$\bar C \equiv  C_*(\widetilde{M}(h'), \tilde E)/i_*(C_*(\tilde X, \tilde A))$
is acyclic;
furthermore, it is free with the evident basis.
Let $\tau(\bar C) \in \Wh_1(K')$ be the Whitehead torsion of $\bar C$
(see \S \chA).
Then, by definition,
$\tau(h) \equiv c_*(\tau(\bar C)) \in \Wh_1(K)$.

By Proposition {\chA.\Alabpc}, %xref
the image of the composite homomorphism:
$$
\pm K \ra K_1(\Z K) @>{\DT}>> HH_1(\Z K)
$$
lies in $HH_1(\Z K)_{C(1)} \subset HH_1(Z K)'$
and thus $\pi'' \circ \DT$ factors through $\Wh_1(K)$
yielding a homomorphism
$\DT''\:\Wh_1(K) \ra HH_1(\Z K)''$.
The isomorphism $h_{\#}\:G \ra K$ induced by $h$ in turn
induces an isomorphism $h_{\#}\:HH_1(\Z G) \ra HH_1(\Z K)$.

\proclaim{Theorem {\chN.\Nlabta} (Torsion Formula)} %xref
For $\gamma \in \Gamma_{(X,A)}$,
$$
 h_{\#}(\tChi''_1(X,A)(\gamma)) - \tChi''_1(Y,B)(h_*(\gamma)) =
(1-h_*(\gamma)) \cdot \DT''(\tau(h)).
$$
\endproclaim

\demo{Proof}
Replacing $h$ with the inclusion of $X$
into the mapping cylinder of $h$,
we may assume without loss of generality that
$h\:(X, A) \hookrightarrow (Y,B)$ is an inclusion of $(X,A)$
into $(Y, B)$ as a subcomplex.
Choose oriented lifts of the cells of $Y$
to the universal cover, $\tilde Y$, of $Y$.
Let $\tilde X = p^{-1}(X)$  and $\tilde A = p^{-1}(A)$
where
$p\:\tilde Y \ra Y$
is the covering projection.
Since $h\:(X, A) \hookrightarrow (Y,B)$ is a homotopy equivalence,
$\tilde X$ is the universal cover of $X$.
Given
\linebreak
$\gamma \in \Gamma_{(X,A)}$,
let $\mu = h_*(\gamma) \in \Gamma_{(Y,B)}$.
We can find,
using the homotopy extension property,
a cellular self homotopy of the identity,
$F^\mu\: (Y,B) \x I \ra (Y,B)$,
representing $\mu$
such that $F^\mu((X,A) \x I) \subset (X,A)$
and the restriction of $F^\mu$ to $(X,A) \x I$
represents $\gamma$.
Let
$(\tilde D_Y^\mu)_* \: C_*(\tilde Y, \tilde B)
\ra C_{* }(\tilde Y,  \tilde B)$
be the chain homotopy determined by the lift of $F^\mu$ and
let $(\tilde D_X^\gamma)_*$ be the restriction of
$(\tilde D_Y^\mu)_*$ to $C_*(\tilde X,  \tilde A)$
(this coincides with the chain homotopy associated
to the lift of the restriction of $F^\mu$ to $(X,A) \x I$).
Let
$
\bar C \equiv C_*(\tilde Y, \tilde B)/h_*(C_*(\tilde X,  \tilde A))$
be the quotient chain complex.
Then $(\tilde D_Y^\mu)_*$ induces a chain homotopy on this
complex which we will denote by $\bar D^\mu_*$.
There is a commutative diagram:
$$
\CD
C_*(\tilde X, \tilde A) @>{}>>  C_{* }(\tilde Y, \tilde B)
@> {} >>  \bar C_* \\
@V{(\tilde D_X^\gamma)_*}VV   @V{(\tilde D_Y^\mu)_*}VV
@V{\bar D^\mu_*}VV \\
C_{* }(\tilde X, \tilde A)   @> {} >>  C_{* }(\tilde Y, \tilde B)
@> {} >>  \bar C_*
\endCD
$$
By \cite{\AJM, Proposition 3.5}, %bibref
we have in $HH_1(\Z K)$:
$$
T_1(\tilde \partial_Y \otimes {\tilde D_Y^\mu}) ~-~
T_1(\tilde \partial_X \otimes {\tilde D_X^\gamma}) =
T_1(\bar \partial \otimes {\bar D^\mu})
\tag{\chN.\Nlabla}
$$
where $\tilde \partial_X$, $\tilde \partial_Y$  and $\bar \del$
are the matrices of the boundary operators of
$C_*(\tilde X, \tilde A)$, $C_{* }(\tilde Y, \tilde B)$ and $C'$ respectively
and $D_X^\gamma$, $D_Y^\mu$ and $\bar D^\mu$ are the matrices
of $(D_X^\gamma)_*$, $(D_Y^\mu)_*$  and $\bar D^\mu_*$ respectively.

Since $h\:(X, A) \hookrightarrow (Y,B)$ is a homotopy equivalence,
$\bar C$ is an acyclic chain complex; furthermore, $\bar C$ is free
with a preferred basis.
In particular, the torsion $\tau(\bar C) \in \Wh_1(K)$ is defined.
By \cite{\AJM, Proposition 3.7}: %bibref
$$
T_1(\bar \partial \otimes {\bar D^\mu}) =
-T_1(V \otimes V^{-1}(1 - \eta_{\#}(\mu)^{-1}))
$$
where $V$ is an invertible matrix representing $\tau(\bar C) = \tau(h)$
(see the proof of
\cite{\AJM, Proposition 7.1}). %bibref
We have
$$
T_1(V \otimes V^{-1}(1 - \eta_{\#}(\mu)^{-1})) =
(1 - \mu) \cdot T_1(V \otimes V^{-1}) = (1 - \mu) \cdot \DT(b)
$$
where $b \in K_1(\Z K)$ is represented by $V$. Also
$$
T_1(\tilde \partial_Y \otimes {\tilde D_Y^\mu}) = \tChi_1(Y,B)(\mu)
\qquad
\text{and}
\qquad
T_1(\tilde \partial_X \otimes {\tilde D_X^\gamma}) =  \tChi_1(X,A)(\gamma).
$$
Substituting into
(\chN.\Nlabla) %xref
and applying $\pi''\:HH_1(\Z K) \ra HH_1(\Z K)''$ yields the
conclusion.\hfill{ }\qed
\enddemo

We say that a homotopy equivalence of pairs $h\:(X,A) \ra (Y,B)$ is
{\it simple} if $\tau(h)=0$
and thus
Theorem {\chN.\Nlabta} %xref
yields:

\proclaim{Theorem {\chN.\Nlabca} (Simple homotopy invariance)} %xref
The derivation 
$$\tChi''_1(X,A)\:\Gamma_{(X,A)} @>{}>> HH_1(\Z G)''$$
is an invariant of the {\it simple} homotopy type of the
pair $(X,A)$, i.e., if $h\:(X,A) \ra (Y,B)$ is a simple homotopy
equivalence then the diagram
$$
\CD
\Gamma_{(X,A)} @>{\tChi''_1(X,A)}>> HH_1(\Z G) \\
@V{h_*}VV         @VV{h_{\#}}V \\
\Gamma_{(Y,B)} @>{\tChi''_1(Y,B)}>> HH_1(\Z K)
\endCD
$$
commutes.\hfill{ }\qed 
\endproclaim

\remark{Remark}
We suspect that
the derivation $\tChi''_1(X,A)$ is {\it not} an invariant
of the homotopy type of the pair $(X,A)$.
This would be the
case if one could find an example
of a homotopy equivalence $h\:(X,A) \ra (Y,B)$ such that the term
$(1-h_*(\gamma)) \cdot \DT''(\tau(h))$ is not
identically zero.
\endremark

A {\it subdivision} $(X',A)$ of a relative CW complex $(X,A)$
is another relative CW structure on the underlying
topological pair $(X,A)$
such that each open cell of $(X',A)$ is contained in some open
cell of $(X,A)$.
The identity  $(X,A) \ra (X',A)$ is a cellular map and a simple
homotopy equivalence. Hence:

\proclaim{Corollary \chN.\Nlabcb} %xref
$\tChi''_1(X,A) = \tChi''_1(X',A)$, i.e.,
$\tChi''_1(X,A)$ is invariant under subdivision.\hfill{ }\qed
\endproclaim

The cohomology classes represented by
the derivations
$\tChi_1(X,A)$,
$\tChi'_1(X,A)$ and
$\tChi''_1(X,A)$
have a stronger invariance property:

\proclaim{Theorem {\chN.\Nlabpf} (Homotopy Invariance)} %xref
The cohomology classes
$\tchi_1(X,A)$,  $\tchi'_1(X,A)$ and
$\tchi''_1(X,A)$ are invariants of the homotopy type of $(X,A)$,
i.e., if $h\:(X,A) \ra (Y,B)$ is a homotopy
equivalence then
$h^*\tchi_1(Y,B)=\tchi_1(X,A)$,
$h^*\tchi'_1(Y,B)=\tchi'_1(X,A)$ and
$h^*\tchi''_1(Y,B)=\tchi''_1(X,A)$.
\endproclaim

\demo{Proof}
The proof given in
\cite{\HEC, Theorem 2.9} %bibref
carries over directly.
Alternatively, 
equation (2.9) decomposes into a pair
equations via the decomposition of 
$HH_1(\Z K)$ 
into $HH_1(\Z K)' \oplus HH_1(\Z K)''$.
The conclusion can then be deduced from
the observation that,
in the proof of
Theorem {\chN.\Nlabta}, %xref
the derivation
$\mu \mapsto
T_1(\bar \partial \otimes {\bar D^\mu}) =(1 - \mu) \cdot \DT(b)$
is inner.\hfill{ }\qed
\enddemo

We can now define
$\tchi_1(N)$,
$\tchi'_1(N)$,
$\tchi''_1(N)$ and
$\tChi''_1(N)$ for a space $N$ which is homeomorphic to
a compact polyhedron:

\definition{Definition \chN.\Nlabdd} %xref
Let $N$ be a path connected topological space which is homeomorphic
to a compact polyhedron.
Let $G= \pi_1(N,x_0)$.
We define a derivation
$$
\tChi''_1(N)\:\Gamma_N \ra HH_1(\Z G)''
\qquad
\text{by}
\qquad
\tChi''_1(N) \equiv  h_{\#} \circ \tChi''_1(|K|) \circ (h^{-1})_*
$$
where $K$ is a finite simplicial complex and
$h\:|K| \ra N$ is a homeomorphism.
We also define cohomology classes by:
$$
\align
\tchi_1(N) &\equiv (h^{-1})^*\tchi_1(|K|) \in H^1(\Gamma_N, HH_1(\Z G)) \\
\tchi'_1(N) &\equiv (h^{-1})^*\tchi'_1(|K|)  \in H^1(\Gamma_N, HH_1(\Z G)') \\
\tchi''_1(N) &\equiv (h^{-1})^*\tchi''_1(|K|)  \in H^1(\Gamma_N, HH_1(\Z G)'').
\endalign
$$
The homomorphisms $\chi_1(N)$, $\chi'_1(N)$ and $\chi''_1(N)$
are similarly defined:
$$
\align
\chi_1(N) &\equiv (h^{-1})^*\chi_1(|K|) \in \Hom(\Gamma_N, H_1(G)) \\
\chi'_1(N) &\equiv (h^{-1})^*\chi'_1(|K|)  \in \Hom(\Gamma_N, H_1(G)) \\
\chi''_1(N) &\equiv (h^{-1})^*\chi''_1(|K|)  \in \Hom(\Gamma_N, H_1(G)).
\endalign
$$
Note that
$\chi_1(N) = \vep_*(\tchi_1(N))$,
$\chi'_1(N) = \vep_*(\tchi'_1(N))$ and
$\chi''_1(N) = \vep_*(\tchi''_1(N))$.
\enddefinition 

By Theorem {\chN.\Nlabca}, %xref
and Chapman's Theorem
that homeomorphisms of compact polyhedra are
simple homotopy equivalences
\cite{Ch}, %bibref
the derivation
$\tChi''_1(N)$ is independent of the choice of $h\:|K| \ra N$.
Also, by Theorem {\chN.\Nlabpf}, %xref
the cohomology classes
$\tchi_1(N)$, $\tchi'_1(N)$ and $\tchi''_1(N)$ are independent
of $h$.

We end this subsection with a technical result,
Proposition {\chN.\Nlabpd}, %xref
which will be used in
\S \chF. %xref

Suppose the topological pair $(Z,C) \subset (X,A)$,
where $C$ is closed in $A$,
is a subcomplex of the finite relative
CW complex $(X,A)$ (i.e., $(Z,C)$ is endowed with the
structure of a relative CW complex such that
each open cell of $(Z,C)$ is an open cell of $(X,A)$).
Assume $Z$ is non-empty, path connected and has a simply connected
universal covering space $\bar Z$.

Define $\bar C = \bar p^{-1}(C)$
where $\bar p\:\bar Z \ra Z$ be the universal covering projection.
Let $\tilde Z = p^{-1}(Z) \subset \tilde X$ and
$\tilde C = p^{-1}(C) \subset \tilde X$.
Let $H=\pi_1(Z,z_0)$ where $z_0 \in Z$ is a chosen basepoint 
for $Z$.
Let $\sigma\:I \ra X$ be a basepath from $z_0$ to $x_0$
(the basepoint of $X$).
The inclusion $i\:Z \hookrightarrow X$ 
and the basepath $\sigma$ determine
a homomorphism $i_{\#}:H \ra G$.
Let $\bar P\: \bar Z \ra \tilde Z$ be a lift of $\bar p$
(i.e., $p \circ \bar P = i \circ \bar p$).
Since $\bar P$ is cellular,
it induces a homomorphism of cellular chain complexes
$\bar P_*\: C_*(\bar Z, \bar C) \ra C_*(\tilde Z, \tilde C)$.
Give $C_*(\bar Z, \bar C) \otimes_{i_{\#}} \Z G$ the right
$G$--module structure $(u \otimes g')\cdot g = u \otimes (g'g)$.
Then $\bar P_*$ induces an isomorphism
$\bar P'_*\: C_*(\bar Z, \bar C)  \otimes_{i_{\#}} \Z G
 \ra C_*(\tilde Z, \tilde C)$
of right $G$--module chain complexes
given by $\bar P'_*(u \otimes g) = \bar P_*(u)g$.
Choose oriented lifts of the cells of $(Z,C)$
to $(\bar Z, \bar C)$.
These lifts and the map $\bar P$ determine
oriented lifts of the cells of
$(Z,C)$ (regarded as a subcomplex of $(X,A)$) to
 $(\tilde Z, \tilde C)$.
Choose oriented lifts of the remaining cells of $(X,A)$ to
$(\tilde X, \tilde A)$.

Suppose $F^\gamma\:(X,A) \x I \ra (X,A)$ is a cellular homotopy
such that  $F^\gamma_0 = F^\gamma_1=\id$ and
$F^\gamma((Z,C) \x I) \subset (Z,C)$.
Let $\bar F^\gamma_{(Z,C)} \:(\bar Z, \bar C) \x I \ra (\bar Z, \bar C)$
be the lift of the restriction of $F^\gamma$ to $(Z,C) \x I$
such that $(\bar F^\gamma_{(Z,C)})_0 = \id$.
Let
$\tilde F^\gamma_{(Z,C)} \: (\tilde Z, \tilde C) \x I \ra (\tilde Z, \tilde C)$
be the restriction of the lift
$\tilde F^\gamma\: (\tilde X, \tilde A) \x I \ra (\tilde X, \tilde A)$
of $F^\gamma$ such that $(\tilde F^\gamma)_0 = \id$.
Let $\bar D^\gamma_*\:C_*(\bar Z, \bar C) \ra C_*(\bar Z, \bar C)$
and $\tilde D^\gamma_*\:C_*(\tilde Z, \tilde C) \ra C_*(\tilde Z, \tilde C)$
be the chain homotopies determined by $\bar F^\gamma_{(Z,C)}$ and
$\tilde F^\gamma_{(Z,C)}$ respectively.
Since $\bar P'_* \circ (\bar D^\gamma_* \otimes \id) = \tilde D^\gamma_*$,
we have $i_{\#}(\bar D^\gamma ) = \tilde D^\gamma$
and
$i_{\#}(\bar \del) = \tilde \del$ (here $\bar D^\gamma$ and
$\bar \del$ are matrices over $\Z H$,
$\tilde D^\gamma$ and
$\tilde \del$ are matrices over $\Z G$,  $\bar \del$ is the
matrix of the boundary operator of $C_*(\bar Z, \bar C)$
and $\tilde \del$ is the
matrix of the boundary operator of $C_*(\tilde  Z, \tilde  C)$).
Thus we obtain the following proposition
which will be useful
in {\S \chF}: %xref

\proclaim{Proposition \chN.\Nlabpd} %xref
For $(Z,C) \subset (X,A)$ as above,
$
i_*(\trace(\bar \del \otimes \bar D^\gamma))=
\trace(\tilde \del \otimes \tilde D^\gamma)
$.\hfill{ }\qed
\endproclaim

\subheading{(C) Geometric Interpretation}

Let $X$ be a compact oriented smooth (or piecewise linear)
manifold.
In cases of interest, $\G (X)$ is non-trivial and so
$\chi(X)=0$ by Gottlieb's Theorem
\cite{Got}. %bibref
Thus, $F^\gamma$ can be perturbed to a map $\bar F^\gamma$ whose
fixed point set,
\hbox{$\Fix(\bar F^\gamma) = \{(x,t) \in
 X \x I ~|~ \bar F^\gamma(x,t)=x\}$}, %WARNING TeXnical trickery
consists of $X \x \{0,1\}$
together with finitely many naturally oriented circles
lying in
\linebreak
$\overset \circ \to X \x (\epsilon, 1-\epsilon)$ for
some $\epsilon > 0$.
Two such circles $V_1$ and $V_2$ lie
{\it in the same fixed point class}
if for some $(x_i,t_i) \in V_i$ there is a path $\nu$ in $X \x I$
from $(x_1,t_1)$ to $(x_2,t_2)$ such that
$(p \circ \nu) (\bar F^\gamma \circ \nu)^{-1}$ is homotopically trivial
where $p\:X \x I \ra X$ denotes projection.
Using $(x_0,0)$ as the basepoint, let $\mu$ be a path from
$(x_0,0)$ to a point $(x,t)$ lying in the circle of fixed points $V$.
We associate with $V$ the conjugacy class $C$ of the ``marker''
$g_C \equiv [(p \circ \mu) (\bar F^\gamma \circ \mu)^{-1}] \in G$;
this conjugacy class is independent of the choice of $\mu$.
Let $\omega$ be the loop based at $(x,t)$ which traverses $V$ once
in the direction of its orientation.
Then $[p \circ (\mu \omega \mu^{-1})] = h$ lies in $Z(g_C)$.
In this way, we also associate
with $V$ an element of
$H_1(Z(g_C)) \cong (Z(g_C))_{\text{ab}}$.
If there are two circles $V_1$ and $V_2$ in the
same fixed point class,
we reach the same centralizer
$Z(g_C)$ from both circles provided the path
used for
$(x_1,t_1) \in V_1$ is $\mu$,
and the path used for $(x_2,t_2) \in V_2$ is $\mu \nu$,
where $\nu$ is such that
$(p \circ \nu) (\bar F^\gamma \circ \nu)^{-1}$ is homotopically trivial.
One treats any (finite) number of
circles similarly.

We form $\hat F^\gamma \: X \x [-\epsilon, 1-\epsilon] \ra X$ by
setting
$$
\hat F^\gamma(x,t) =
\cases
\bar F^\gamma(x,t) &\text{if $t \in [0,1-\epsilon]$} \\
\bar F^\gamma(x,1+t) &\text{if $t \in [-\epsilon, 0]$.}
\endcases
$$
Perturbing $\hat F^\gamma$ near $X \x \{0\}$,
we may assume that $\Fix(\hat F^\gamma)$ is a union of naturally
oriented circles: the previous circles and some new ones
near $X \x \{0\}$.
With each {\it new} circle we associate the trivial conjugacy class
and an element $[p \circ (\mu\omega\mu^{-1})]$, as before,
lying in $H_1(G)\equiv H_1(Z(1))$.

In this way we obtain from $F^\gamma$ a slightly altered map
$\hat F^\gamma$, and we associate with $\Fix(\hat F^\gamma)$
an element of $\bigoplus_{C \in G_1} H_1(Z(g_C)) \cong HH_1(\Z G)$;
for details see
\cite{\TTTF; \S 1},  %bibref
\cite{\HEC; \S 1},   %bibref
and
\cite{\AJM; \S 6}.   %bibref

If we only consider (oriented) circles whose associated conjugacy
classes lie $\G(X)$ [respectively $G_1 - \G (X)$] we get
geometric interpretations of $\tChi'_1(X)(\gamma)$
[respectively $\tChi''_1(X)(\gamma)$].
By interpreting these oriented circles as homology classes in $X \x I$,
using $H_1(X \x I) \cong H_1(X) \cong H_1(G)$, we get
$\chi'_1(X)(\gamma)$ and $\chi''_1(X)(\gamma)$.

%%%%%%%%%%%%%%%%%%%%%%%%%%
%%%   sec03.tex           
%%%%%%%%%%%%%%%%%%%%%%%%%%

\heading
\S \chF. Filtered cell complexes
\endheading

The appropriate notion of a ``$G$--CW complex''
\cite{I$_1$}, %bibref
where $G$ is a compact Lie group,
requires us to consider spaces obtained by attaching
cells which a priori are more general than CW-complexes.
In this paper
we will only be concerned with the cases $G=S^1$, the circle group,
and $G=T^2\equiv S^1 \x S^1$, the $2$--torus.
The main result of this section is
Theorem {\chF.\Flabta} %xref
which will be used in {\S \chS} %xref
(see Theorem {\chS.\Slabta}) %xref
to establish the connection between the topological
invariants of {\S \chN} %xref
and the $S^1$-Euler characteristic.

\definition{Definition \chF.\Flabda}    %xref
A {\it filtered cell complex} is a compact Hausdorff space
$X$ together with a filtration 
$\emptyset \equiv X_{-1} \subset X_0 \subset \cdots \subset X_N = X$
by closed subsets such that for each
\linebreak
$n=0,\ldots N$ the pair $(X_n, X_{n-1})$ is a finite relative CW complex.
\enddefinition

Note that $X_0$ is a finite CW complex.

A filtered cell complex is a compact ANR
\cite{Wh}; %bibref
however, as we haven't fully developed $1$--parameter fixed point theory
for general compact ANR's, it is convenient to deal with a subclass
of such spaces:

\definition{Definition \chF.\Flabdd}    %xref
A filtered cell complex $X$ is {\it polyhedral} if $X_0$ is
a compact polyhedron and all attaching maps in the relative
CW complexes $(X_n,X_{n-1})$ are piecewise linear.
\enddefinition
In particular, $X$ is itself a polyhedron.

Given a filtered cell complex $X$
with filtration $\F = \{X_i\}$,
let $(X; \F)^{(X; \F)}$ denote the subspace of $X^X$ consisting
of filtration preserving maps, i.e., maps $f\:X \ra X$  such that
$f(X_i) \subset X_i$ for all $i$.
Define $\Gamma_{(X; \F)} \equiv \pi_1((X; \F)^{(X; \F)},\id)$.
For $j \geq i$, restriction to $X_j$ yields a map
$(X; \F)^{(X; \F)} \ra (X_j, X_i)^{(X_j, X_i)}$ which in turn
induces a homomorphism  $\Gamma_{(X; \F)} \ra \Gamma_{(X_j, X_i)}$.
We abuse notation: given
$\gamma \in \Gamma_{(X; \F)}$ we also denote by $\gamma$ the image
of $\gamma$ in $\Gamma_{(X_j, X_i)}$.

By an inductive application of the relative
Cellular Approximation Theorem
and the homotopy extension property,
any $\gamma \in \Gamma_{(X; \F)}$ can be represented by a
homotopy
\linebreak
$F^\gamma\:X \x I \ra X$ such that:
\roster
\item
$F^\gamma_0=F^\gamma_1=\id$.
\item
$F^\gamma$ is filtration preserving,
i.e.,  $F^\gamma(X_n \x I) \subset X_n$ for
all $n$.

\item
For each $n$, the map of pairs
$F^{\gamma, n, n-1} \:(X_n, X_{n-1}) \x I \ra (X_n, X_{n-1})$
given by the  restriction of $F^\gamma$ to $X_n$ is cellular.
\endroster

Suppose $X$ is path connected.
Let $\tilde X_n$ be the inverse image of $X_n$ in  $\tilde X$,
the universal cover of $X$.
Let $G \equiv \pi_1(X,x_0)$.
For each $n$ and for each cell $e$ of $(X_n, X_{n-1})$,
choose an oriented lift of $e$ to $\tilde X$.
Let  $\tilde F^\gamma\:\tilde X \x I \ra \tilde X$ be the
unique lift of $F^\gamma$ such that $\tilde F^\gamma_0 = \id$ 
(recall that $\tilde F^\gamma_1$ is the covering transformation
corresponding to $\eta_{\#}(\gamma) \in G$).
Since
$\tilde F^{\gamma, n, n-1}\:
(\tilde X_n, \tilde X_{n-1}) \x I \ra (\tilde X_n, \tilde X_{n-1})$
given by the  restriction of $\tilde F^\gamma$ to $\tilde X_n$ is cellular,
we obtain chain homotopies
(retaining the Sign Convention 
of {\S \chN}): %xref
$$
\tilde D^{\gamma, n, n-1}_*\: C_*(\tilde X_n, \tilde X_{n-1})
\ra  C_*(\tilde X_n, \tilde X_{n-1})
\qquad n \geq 0.
$$
Let $\tilde \del^{n,n-1}_*$
denote the boundary operator in $C_*(\tilde X_n, \tilde X_{n-1})$.

\definition{Definition \chF.\Flabdc} %xref
Define the function $\tChi_1(X;\F)\:\Gamma_{(X; \F)} \ra HH_1(\Z G)$ by
$$
\tChi_1(X; \F)(\gamma) \equiv
\sum_{n \geq 0}  T_1(\tilde \del^{n,n-1} \otimes \tilde  D^{\gamma, n, n-1} ).
$$
\enddefinition

It is straightforward to verify that,
for a given choice of lifts of cells,
$\tChi_1(X; \F)(\gamma)$ does not depend on the the choice of
$F^\gamma$ representing $\gamma$ and that 
the function $\tChi_1(X; \F)$ is a derivation.
As in {\S \chN}, %xref
we define derivations:
$$
\tChi'_1(X; \F)\:\Gamma_{(X; \F)} \ra HH_1(\Z G)'
\qquad \text{and} \qquad
\tChi''_1(X; \F)\:\Gamma_{(X; \F)} \ra HH_1(\Z G)''
$$
by
$\tChi'_1(X; \F) \equiv \pi' \circ  \tChi_1(X; \F)$ and
$\tChi''_1(X; \F) \equiv \pi'' \circ \tChi_1(X,A)$.
Observe that
\linebreak
$\tChi_1(X; \F) =  (\tChi'_1(X; \F), \tChi''_1(X; \F))$.

The term $T_1(\tilde \del^{n,n-1} \otimes \tilde  D^{\gamma, n, n-1} )$
appearing in
Definition  {\chF.\Flabdc} %xref
can be identified as follows.
Let $X_{n,1}, \ldots, X_{n,k_n}$ be the path components of $X_n$.
For $j=1,\dots,k_n$,
choose basepoints $x_{n,j} \in X_{n,j}$ and
basepaths $\sigma_{n,j}\:I \ra X$ from $x_{n,j}$ to $x_0$ (the
basepoint for $X$).
Let $H_{n,j} \equiv \pi_1(X_{n,j}, x_{n,j})$ and let
$i_j\:H_{n,j} \ra G$ the homomorphism determined by the
inclusion $X_{n,j} \subset X$ and the basepath $\sigma_{n,j}$.
By
Proposition {\chN.\Nlabpd}, %xref
we have
(after choosing cell lifts in the manner prescribed
in the discussion preceding 
that proposition): %xref

\proclaim{Proposition \chF.\Flabpa} %xref
For $n \geq 0$:
$$
T_1(\tilde \del^{n,n-1} \otimes \tilde  D^{\gamma, n, n-1} )
=
\sum^{k_n}_{j=1}
(i_j)_*(\tChi_1(X_{n,j}, X_{n,j} \cap X_{n-1})(\gamma)).\hfill{ }\qed
$$
\endproclaim

\proclaim{Proposition \chF.\Flabpb} %xref
Suppose $X$ is a finite connected CW complex and $\F =\{X_i\}$ is
a  filtration of $X$ by subcomplexes.
Then for $\gamma \in \Gamma_{(X; \F)}$,
$\tChi_1(X)(\gamma) = \tChi_1(X; \F)(\gamma)$.
\endproclaim

\demo{Proof}
Clearly,
$X$ together with $\F$ is a filtered cell complex.
Recall that $\tilde X_n = p^{-1}(X_n)$ where $p\:\tilde X \ra X$
is the universal covering projection.
Choose oriented lifts of each cell of $X$ to $\tilde X$.
Let $\tilde X^k_n$ be the $k$--skeleton of $\tilde X_n$.
The Mayer-Vietoris sequence yields a short exact sequence:
$$
0 \ra H_k(\tilde X_{n-1}^k, \tilde X_{n-1}^{k-1}) \ra
H_k(\tilde X_n^k, \tilde X_n^{k-1}) \ra
H_k(\tilde X_n^k \cup \tilde X_{n-1}, \tilde X_n^{k-1}\cup \tilde X_{n-1}) 
\ra 0
$$
and thus for each $n \geq 0$ we have an exact sequence of cellular
chain complexes:
$$
0 \ra C_*(\tilde X_{n-1}) \ra C_*(\tilde X_n)
\ra C_*(\tilde X_n, \tilde X_{n-1}) \ra 0.
$$
Represent $\gamma \in \Gamma_{(X; \F)}$ by a cellular map
$F^\gamma \: X \x I \ra X$ such that
$F^\gamma_0 = F^\gamma_1 = \id$ and 
$F^\gamma(X_n \x I) \subset X_n$ for all $n$.
Let $\tilde F^{\gamma, n}\:\tilde X_n \x I \ra \tilde X_n$
denote the restriction of
$\tilde F^\gamma\: \tilde X \x I \ra \tilde X$
to $\tilde X_n$ and let
$
\tilde D^{\gamma, n}_*\: C_*(\tilde X_n) \ra  C_*(\tilde X_n).
$
be the corresponding chain homotopy.
Let $\tilde \del^n_*$
denote the boundary operator in $C_*(\tilde X_n)$.
There is a commutative diagram:
$$
\CD
C_*(\tilde X_{n-1}) @>{}>>  C_*(\tilde X_n)
@> {} >>  C_*(\tilde X_n, \tilde X_{n-1}) \\
@V{\tilde D^{\gamma, n-1}_*}VV   @V{\tilde D^{\gamma, n}_*}VV
@V{\tilde D^{\gamma, n, n-1}_*}VV \\
C_*(\tilde X_{n-1}) @>{}>>  C_*(\tilde X_n)
@> {} >>  C_*(\tilde X_n, \tilde X_{n-1}).
\endCD
$$
By \cite{\AJM, Proposition 3.5}, %bibref
$$
T_1(\tilde \del^{n, n-1} \otimes \tilde D^{\gamma, n, n-1}) =
T_1(\tilde \del^n \otimes \tilde D^{\gamma, n}) -
T_1(\tilde \del^{n-1} \otimes \tilde D^{\gamma, n-1}).  
$$
Substituting into
Definition {\chF.\Flabdc}, %xref
$$
\tChi_1(X; \F)(\gamma) =
\sum_{n \geq 0}  T_1(\tilde \del^{n,n-1} \otimes \tilde  D^{\gamma, n, n-1}) =
\sum_{n \geq 0} (
T_1(\tilde \del^n \otimes \tilde D^{\gamma, n}) -
T_1(\tilde \del^{n-1} \otimes \tilde D^{\gamma, n-1})).
$$
The last sum collapses to
$$
T_1(\tilde \del^N \otimes \tilde D^{\gamma, N}) =
T_1(\tilde \del \otimes \tilde D^{\gamma}) =
\tChi_1(X)(\gamma)
$$
where $X_N=X$.\hfill{ }\qed
\enddemo

The following theorem is the main technical ingredient in
the proof of Theorem {\chS.\Slabta} %xref
where it will be applied to an $S^1$--CW complex which is
filtered by its $S^1$--skeleta.

\proclaim{Theorem \chF.\Flabta} %xref
Let $X$ be a path connected polyhedral filtered cell
complex with filtration $\F$.
Then
$
\tChi''_1(X) \circ i = \tChi''_1(X; \F) 
$
and
$
i^*\tchi'_1(X) = [\tChi'_1(X; \F)] 
$
where $i\:\Gamma_{(X; \F)} \ra \Gamma_X$ is the natural homomorphism.
\endproclaim

\demo{Proof}
Since the filtered cell complex $(X; \F)$ is polyhedral
(see Definition \chF.\Flabdd), %xref
$X$ is a compact polyhedron and there is a triangulation
$X'$ of $X$ and a filtration $\F' \equiv \{X'_n\}$ of $X'$ by subcomplexes
such that for each $n$, $(X'_n, X_{n-1})$ is a subdivision
of the relative CW complex $(X_n, X_{n-1})$.
By Proposition {\chF.\Flabpb}, %xref
$\tChi_1(X') \circ i = \tChi_1(X'; \F')$.

With notation as in
Proposition {\chF.\Flabpa}, for each $n \geq 0$
let $X_{n,j}$, $j=1,\ldots, k_n$ be the set of
path components of $X_n$; let $X'_{n,j}$ be the
triangulation of $X_{n,j}$ determined by $X'_n$.
The proof of
Theorem {\chN.\Nlabta} %xref
shows that for $\gamma \in \Gamma_{(X; \F)}$:
$$
\tChi_1(X_{n,j}, X_{n,j} \cap X_{n-1})(\gamma)
~-~
\tChi_1(X'_{n,j}, X_{n,j} \cap X_{n-1})(\gamma)
= (1 - \gamma) \cdot \DT(b_{n,j})
$$
where $b_{n,j} \in K_1(\Z H_{n,j})$ is a representative of
the torsion of the identity map
\linebreak
$\id_{n,j}\:(X_{n,j}, X_{n,j} \cap X_{n-1}) \ra
(X'_{n,j}, X_{n,j} \cap X_{n-1})$.
Each $\id_{n,j}$ is a simple homotopy equivalence and thus
$\DT(b_{n,j}) \in HH_1(\Z H_{n,j})_{C(1)}$.
It follows that 
$(1 - \gamma) \cdot (i_j)_*\DT(b_{n,j}) \in HH_1(\Z G)'$.
By Proposition {\chF.\Flabpa},
$$
\tChi_1(X; \F)(\gamma) - \tChi_1(X'; \F')(\gamma)
=  (1 - \gamma) \cdot
\left( \sum_{n,j} (i_j)_*\DT(b_{n,j}) \right).
$$
    From this we deduce that
$\tChi''_1(X; \F)= \tChi''_1(X'; \F')$
and that
$[\tChi'_1(X; \F)]= [\tChi'_1(X'; \F')]$.
Now by
Definition {\chN.\Nlabdd}, %xref
$\tChi''_1(X)  = \tChi''_1(X')$ and
$\tchi'_1(X)   = \tchi'_1(X')$.\hfill{ }\qed
\enddemo

%%%%%%%%%%%%%%%%%%%%%%%%%%
%%%   sec04.tex           
%%%%%%%%%%%%%%%%%%%%%%%%%%

\heading
\S \chS. The $S^1$--Euler characteristic of an $S^1$--CW complex
\endheading

In this section we introduce
the ``$S^1$--Euler characteristic''
(Definition \chS.\Slabdc) %xref
of a finite $S^1$--CW complex $X$ and show in
Proposition {\chS.\Slabpc} %ref
that it coincides with
$\tChi_1(X; \F)(\gamma)$
where $\gamma \in \Gamma_X$ is the element defined
by the $S^1$--action and $\F$ is the filtration of
$X$ by
\linebreak
$S^1$--skeleta.
This, together with
Theorem {\chF.\Flabta}, %xref
yields
Theorem {\chS.\Slabta}  %xref
which establishes the precise relationship between
the $S^1$--Euler characteristic and the topological invariants
$\tchi'_1(X)$ and $\tChi''_1(X)$ of 
\S \chN. %xref
In addition,
we are led to a concise formula
(Theorem \chS.\Slabtb) %xref
for $\chi_1(X)(\gamma)$.
We give two applications of these theorems to manifolds.
Theorem {\chS.\Slabtf} %xref
asserts that if $X$ is a closed even dimensional smooth $S^1$--manifold
then $\chi_1(X)(\gamma)=0$.
Theorem {\chS.\Slabte}  %xref
gives a formula
for the Poincar\'e dual of the
Euler class of the normal bundle to the flow defined
by a smooth $S^1$--action
without fixed points on a smooth closed oriented manifold.
Theorem {\chS.\Slabtg} %xref
establishes the formulas (0.4) %xref
and (0.5) %xref
for $3$--dimensional Seifert fibered spaces.

Let $S^1$ denote the set of complex numbers of unit modulus
regarded as a compact Lie group.
An {\it $S^1$--space} is a topological space $X$ together with
a continuous left action, $S^1 \x X @>{\alpha}>> X$, of $S^1$ on $X$.

We recall from \cite{I$_1$} %bibref
the notion of an ``$S^1$--CW complex''.

\definition{Definition \chS.\Slabda}    %xref
Let $X$ be a Hausdorff $S^1$--space, $A$ a closed $S^1$--subset
of $X$  (i.e.,
\linebreak
$\alpha(S^1 \x A) \subset A$) 
and $n$ a non-negative integer.
We say
{\it $X$ is obtained from $A$ by attaching $S^1$--$n$--cells}
if there is a collection
$\{c^n_j ~|~ j \in J\}$ of closed $S^1$--subsets of $X$ such that:
\roster
\item
$X = A \cup \bigcup_{j \in J} c^n_j$, and $X$ has the topology
coherent with $A$ and $\{c^n_j ~|~ j \in J\}$.

\item
For $i \neq j$,
$(c^n_i - A) \cap (c^n_j - A) = \emptyset$.

\item
For each $j \in J$ there is a closed subgroup $H_j$ of $S^1$ 
(note that $H_j$ is finite or $H_j=S^1$)
and
an $S^1$--map
$f_j\:(S^1/H_j \x D^n, S^1/H_j \x S^{n-1}) \ra (c^n_j, c^n_j \cap A)$
such that $f_j(S^1/H_j \x D^n) = c^n_j$ and $f_j$ maps
$S^1/H_j \x D^n  ~-~ S^1/H_j \x S^{n-1}$ homeomorphically onto
$c^n_j - A$.
\endroster
Each $c^n_j$ is called an {\it $S^1$--$n$--cell}.
The restriction  $f_j|\:S^1/H_j \x S^{n-1} \ra A$ of $f_j$ is
called the {\it attaching map} of the $S^1$--$n$--cell $c^n$.
\enddefinition

\definition{Definition \chS.\Slabdb}    %xref
An
{\it $S^1$--relative CW complex}
$(X,A)$ consists of a Hausdorff
$S^1$--space $X$, a closed $S^1$--subset $A$ of $X$ and an increasing
filtration of $X$ by closed $S^1$--subsets
$(X,A)_k$, $k=0,1,\ldots$, such that:
\roster
\item
$(X,A)_0$ is obtained from $A$ by attaching $S^1$--$0$--cells, and
for $k \geq 1$, $(X,A)_k$ is obtained from $(X,A)_{k-1}$ by attaching
$S^1$--$k$--cells.
\item
$X = \bigcup_{k \geq 0} (X,A)_k$, and $X$ has the topology
coherent with $\{(X,A)_k ~|~ k \geq 0\}$.
\endroster
\enddefinition
The closed $S^1$--subset $(X,A)_k$ is called the
{\it $S^1$--$k$--skeleton of $(X,A)$}.
If $A = \emptyset$ we call $X$ an 
{\it $S^1$--CW complex}
and denote
its $S^1$--$k$--skeleton by $X_k$.
For convenience, we adopt the convention $X_{-1} = \emptyset$.
We say $X$ is {\it finite} if $X$ has finitely many $S^1$--cells.

\proclaim{Proposition \chS.\Slabpd}  %xref
Let $X$ be a finite $S^1$--CW complex.
Then $(X, \{X_k\})$ is a filtered cell complex
(Definition \chF.\Flabda) %xref
in which
the relative CW complex structure on $(X_k,X_{k-1})$ is
related to the $S^1$--$k$--cells $c^k_j$ as follows:
when $H_j$ is finite, $c^k_j$ is covered by one
$(k+1)$--cell $d^{k+1}_j$ and one $k$--cell $e^k_j$ of
$(X_k,X_{k-1})$; when $H_j = S^1$,
$c^k_j$ is covered by one
$k$--cell $e^k_j$.
The images of the cells $\{e^k_j\}$ decompose the
quotient space $X/S^1$ as a CW complex.
\endproclaim

\demo{Proof}
Let $Y$ be a space and $f\:S^1 \x S^{k-1} \ra Y$ a map and
$f_| \:\{z_0\} \x S^{k-1} \ra Y$
\linebreak
its restriction where
$z_0 \in S^1$ is a basepoint.
Let $Z = (S^1 \x D^k) \cup_f Y$  and $Z' = (\{z_0\} \x  D^k) \cup_{f_|} Y$
so that $Z' \subset Z$.
Define $F\:\del(I \x D^k) \ra Z'$ by
$F(i,x)=(z_0,x)$ if $i=0,1$ and $x \in D^k$, and
$F(t,y) = f(\exp(2\pi \imath t), y)$ if $y \in S^{k-1}$ and $t \in I$.
Then there is a natural homeomorphism $Z \cong (I \x D^k) \cup_F Z'$.
Thus $(Z,Y)$ is a relative CW complex, obtained by attaching a
$k$--cell $e$ and then a $(k+1)$--cell $d$.
Our claim for the case where $H_j$ is finite follows.
The case $H_j=S^1$ is similar.
Although $X$ is thus obtained by attaching cells, it might fail
to be a CW complex only because the cells $e^k_j$ might
be attached to the $k$--skeleton rather than the $(k-1)$--skeleton;
however, this difficulty disappears in the quotient space.\hfill{ }\qed
\enddemo

Theorem 7.1 and Proposition 6.1 of \cite{I$_2$}  %bibref
give:

\proclaim{Proposition {\chS.\Slabtd}} %xref
Any smooth $S^1$--manifold $M$ (with or without boundary)
has an
\linebreak
``$S^1$--triangulation'' and thus an
$S^1$--CW complex structure (which is finite if $M$ is compact)
such that the attaching maps of the $S^1$--cells are
piecewise linear.
In particular,
an $S^1$--triangulation gives rise to a polyhedral
filtered cell complex
(Definition \chF.\Flabdd).\hfill{ }\qed %xref
\endproclaim

Let $X$ be a connected finite $S^1$--CW complex with base point $x_0$,
and let
\linebreak
$G = \pi_1(X,x_0)$.
The
action $\alpha \:S^1 \x X \rightarrow X$ is adjoint to
a loop in $X^X$ based at $\id_X$, and hence defines 
a {\it fundamental element}
$\gamma \in \Gamma_X$.
Define $F^\gamma \: X \x I \to X$ by
$F^\gamma(x,t) = \alpha(e^{2\pi \imath t}, x)$.
Let the $S^1$--$n$--cells be $\{c^n_j \mid j\in J_n\}$.  

With each $c^n_j$ we associate the {\it primitive loop}
$\omega_j \: S^1 \to X$:
when $H_j$ is finite, the formula is
$\omega_j(e^{2\pi \imath t}) = f_j(e^{2\pi \imath t/|H_j|}, 0)$,
an embedding of $S^1$ which goes around the $S^1$--orbit of
$\omega_j(1) \equiv f_j(1,0)$ once.
When $H_j = S^1$, 
$\omega_j$ is the constant loop
at $f_j(S^1/S^1, 0)$.  
It is convenient to set $|H_j| = 0$ when $H_j = S^1$.
Then the loops $\omega_j^{|H_j|}$ are freely homotopic
for all $n$ and all $j \in J_n$.  

The choice of a lift $\tilde e^n_j$ or $\tilde d^{n+1}_j$ implies the choice
of a base path $\sigma_j$ (up to homotopy) for the cell $e^n_j$ or
$d^{n+1}_j$, with $\sigma_j(0) = x_0$.
When $H_j$ is finite, we pick
abutting lifts $\tilde e^n_j$ and $\tilde d^{n+1}_j$ with the same base
path $\sigma_j$; more precisely
(see the Sign Convention in \S \chN)
we arrange:
$$
\tilde{\partial}^{n,n-1}_{n+1}(\tilde d^{n+1}_j) = (-1)^n\tilde
e^n_j ~(g^{-1}_{j,n} - 1).
\tag{\chS.\Slabla}  %xref
$$
where $g_{j,n} \in G$ is represented by the loop
$\sigma_j \omega_j \sigma^{-1}_j$.
We have $g_{j,n}^{|H_j|} =\eta_{\#}(\gamma)$.
In particular,
if any $H_j=S^1$ then $\eta_{\#}(\gamma)=1$.

\definition{Definition {\chS.\Slabdc}  %xref
$(S^1$--Euler characteristic)} 
$\tchi_{S^1}(X) \in HH_1(\Z G)$ is the element represented by the
Hochschild $1$--cycle
$$
\sum_{n\geq 0} (-1)^{n+1}\sum_{j\in J_n} \sum^{|H_j|}_{i=1} g_{j,n}
\otimes g^{-1}_{j,n}g^{-i}_{j,n}.
$$
\enddefinition

\remark{Remark}
$\tchi_{S^1}(X)$ is independent of the choice of the basepaths $\sigma_j$.
Another choice of basepaths gives rise to $g'_{j,n} \in G$ related to
$g_{j,n}$ by $g'_{j,n}= h_{j,n} g_{j,n} h^{-1}_{j,n}$ for some
$h_{j,n} \in G$.
Then the corresponding cycles in
Definition {\chS.\Slabdc} %xref
are homologous because
$$
g'_{j,n} \otimes (g'_{j,n})^{-1-i} - g_{j,n} \otimes g^{-1-i}_{j,n}
~=~
d \left(
h^{-1}_{j,n} \otimes h_{j,n} g_{j,n} \otimes g^{-1-i}_{j,n} -
h_{j,n}  g_{j,n} \otimes h^{-1}_{j,n}
\otimes (h_{j,n} g_{j,n} h^{-1}_{j,n})^{-1-i}
\right).
$$
\endremark

By Proposition {\chS.\Slabpd}, %xref
we may regard $X$ as a finite filtered cell complex $(X,\F)$,
where $\F=\{X_n\}$ is the filtration by $S^1$--$n$--skeleta.
The cellular chain complex $C_*(\tilde X_n,\tilde X_{n-1})$
is zero in degrees other than $n$ and $n+1$:
$$
0 \ra
C_{n+1}(\tilde X_n,\tilde X_{n-1})
@>{\tilde\partial^{n,n-1}_{n+1}}>>
C_n(\tilde X_n,\tilde X_{n-1})
\ra 0.
$$
The chain homotopy $\tilde D^{\gamma,n,n-1}_*$ is given by
$$
\tilde D^{\gamma,n,n-1}_n(\tilde e^n_j) =
(-1)^{n+1}\tilde d^{n+1}_j \left(1 + g^{-1}_{j,n} +
\cdots + g^{-(|H_j| - 1)}_{j,n} \right)
\tag{\chS.\Slablb} %xref
$$
(interpreted as $0$ when $|H_j| = 0$).
By our conventions,
$\tilde D^{\gamma,n,n-1} = (-1)^{n+1}\tilde D^{\gamma,n,n-1}_n$
and
$\tilde\partial^{n,n-1} = \tilde\partial^{n,n-1}_{n+1}$.  

\proclaim{Proposition \chS.\Slabpc} %xref
$\tchi_{S^1}(X) = \tChi_1(X;\F)(\gamma)$
\endproclaim

\demo{Proof}
    From Definition {\chF.\Flabdc} %xref
and the expressions (\chS.\Slabla) %xref
and (\chS.\Slablb), %xref
we see that $\tChi_1(X;\F)(\gamma)$ is represented by 
$$
\sum_{j\geq 0}(-1)^n(g^{-1}_{j,n}-1)\otimes
\bigg(
\sum_{i=0}^{|H_j|-1}g^{-i}_{j,n}
\bigg).
$$ 
By standard Hochschild chain identities
(see \cite{\AJM; Lemma 6.13}), %bibref
a term of the form $1 \otimes g$ is
homologically trivial, and a term $g^{-1}\otimes gg^{-i}$ is homologous to
$-g\otimes g^{-1}g^{-i}$.\hfill{ }\qed
\enddemo

Applying
Theorem {\chF.\Flabta}, %xref
to the right side of the
equality $\tchi_{S^1}(X) = \tChi_1(X;\F)(\gamma)$
given by
Proposition {\chS.\Slabpc}, %xref
we obtain.

\proclaim{Theorem \chS.\Slabta}  %xref
If $X$ is a connected finite $S^1$--CW complex
which is polyhedral as a filtered cell complex and
$\gamma \in \Gamma_X$ is the fundamental element then:
\roster
\item
$\tChi''_1(X)(\gamma) =
\pi''(\tchi_{S^1}(X))$.

\item
Let
$i\:\langle \gamma \rangle \hookrightarrow \Gamma_X$
be the inclusion of the cyclic subgroup generated by
$\gamma$.
Then the cohomology class
$i^*\tchi'_1(X) \in H^1(\langle \gamma \rangle, HH_1(\Z G)')$
is represented by the derivation
$d_X\:\langle \gamma \rangle \ra  HH_1(\Z G)'$,
$d_X \equiv \tChi_1(X;\F) \circ i$,
which takes  $\gamma$
to $\pi'(\tchi_{S^1}(X))$.\hfill{ }\qed
\endroster
\endproclaim

\remark{Remark}
Definition {\chS.\Slabdc}  %xref
has both a ``combinatorial'' and a geometric motivation.
Note that $\tchi_{S^1}(X)$ is defined in terms of the
$S^1$--CW structure of $X$ (and so may be viewed as
combinatorial); the defining formula is motivated by
Proposition {\chS.\Slabpc}. %xref
Theorem {\chF.\Flabta} %xref
implies that $\tchi_{S^1}(X)$ is essentially the same
as $\tChi_1(X')(\gamma)$  where $X'$ is a CW subdivision
of $X$. If $X'$ is a PL manifold then $\tChi_1(X')(\gamma)$
has a natural geometric interpretation in terms of fixed
point theory and transversality as described in
\S \chN(C). %xref
\endremark
\smallskip

\flushpar
{\bf Notation.}
We write $\{g\} \in H_1(G) \equiv G_{ab}$ for $A(g)$, the image of
$g\in G$ under abelianization.
\smallskip

The following theorem generalizes \cite{\HEC; Theorem 4.5} %bibref
where the case of a free $S^1$--action is treated.

\proclaim{Theorem \chS.\Slabtb}   %xref
Let $X$ be as in
Theorem {\chS.\Slabta},  %xref
If the $S^1$--action has a fixed point then
$\chi_1(X)(\gamma) = 0$.  If there are no fixed points then
$\chi_1(X)(\gamma) = -\chi(X/S^1)\{\eta_{\#}(\gamma)\}$.
\endproclaim

\demo{Proof}
We have $|H_j|\{g_{j,n}\} = \{\eta_{\#}(\gamma)\}$.
By Propositions
{\chA.\Alabpa}, %xref
{\chA.\Alabpb}  %xref
and
{\chS.\Slabpc},  %xref
$$
\align
\chi_1(X)(\gamma) &= \sum_{n\geq 0}(-1)^{n+1}\sum_{j\in J_n}
|H_j|\{g_{j,n}\} \\
&= \sum_{n\geq 0} (-1)^{n+1}\sum_{j\in J_n}\{\eta_{\#}(\gamma)\}
\endalign
$$
If there is a fixed point then $\eta_{\#}(\gamma) = 1$, so the right hand
side is 0.  If there is no fixed point, then the obvious bijection between
the $S^1$--$n$--cells of $X$ and the $n$--cells of $X/S^1$
(Proposition \chS.\Slabpd) %xref
makes the right hand
side $-\chi(X/S^1)\{\eta_{\#}(\gamma)\}$.\hfill{ }\qed
\enddemo

Recall that there is a natural homomorphism $\vep_*\: HH_1(\Z G) \ra H_1(G)$
which takes a Hochschild homology class represented by
$\sum_i n_i g_i \otimes h_i$ to $\sum_i n_i \{g_i\}$
(see {\S \chA}). %xref
The proof of Theorem {\chS.\Slabtb}   %xref
shows that:

\proclaim{Addendum} $\vep_*(\tchi_{S^1}(X)) = \chi_1(X)(\gamma)$.
\endproclaim

We give two application to manifolds.

Our first application generalizes the classical theorem that
the Euler characteristic of a closed odd dimensional manifold
(whether orientable or not) is zero.

\proclaim{Theorem \chS.\Slabtf}  %xref
Let $X$ be a closed even dimensional smooth manifold and
suppose $\gamma  \in \Gamma_X$ is represented by a
smooth $S^1$--action on $X$.
Then $\chi_1(X)(\gamma)=0$.
\endproclaim

\demo{Proof}
By Proposition {\chS.\Slabtd} %xref
there exists an $S^1$--CW complex structure on the
\linebreak
$S^1$--manifold $X$.
If the $S^1$--action on $X$ has a fixed point then
by the first part of
\linebreak
Theorem {\chS.\Slabtb} %xref
we have
$\chi_1(X)(\gamma)=0$.
If the $S^1$--action has no fixed points then the
quotient $X/S^1$ is a rational homology
manifold of odd dimension and so
$\chi(X/S^1) = 0$.
Hence,
by Theorem {\chS.\Slabtb}, %xref
$\chi_1(X)(\gamma) = -\chi(X/S^1)\{\eta_{\#}(\gamma)\} = 0$.
\enddemo

\remark{Remark}
Theorem {\chS.\Slabtf}  %xref
generalizes
\cite{GNO, Corollary 7.5} %bibref
where, using different techniques,
the same conclusion (with rational coefficients) was obtained
under the additional assumption that $X$ is symplectic and
the $S^1$--action is Hamiltonian.
\endremark

Our second application
(Theorem {\chS.\Slabte}) %xref
gives a formula for the Poincar\'e dual of the
Euler class of the normal bundle to the flow defined
by a smooth $S^1$--action
without fixed points on a closed oriented manifold.
Let $X$ be a connected smooth closed oriented
$S^1$--manifold partitioned as a finite
$S^1$--CW complex which is polyhedral as a filtered
cell complex
(see Proposition \chS.\Slabtd). %xref
Assume that the action has no fixed points,
so that each $H_j$ is finite.
Let $\lambda$ be the real line bundle over $X$
consisting of tangent vectors which are tangent
to the $S^1$--orbits.
The $S^1$--action determines an orientation on $\lambda$.
Let $\nu = T_X/\lambda$ where $T_X$ is the tangent
bundle of $X$.
Given a Riemannian metric on $X$, $\nu$ is identified
with the oriented normal bundle to the flow defined by
the $S^1$--action.
Recall from
\cite{\HEC}  %bibref
that (in view of our use of Dold's sign conventions
for cap products)
Poincar\'e duality is given by
$\PD_X(u) = (-1)^{i(m-i)} u \cap [X]$ where
$u \in H^i(X)$ and $[X] \in H_m(X)$ is the fundamental
class of $X$.
We compute the Poincar\'e dual of the Euler class
$\euler(\nu) \in H^{m-1}(X; \Z)$:

\proclaim{Theorem \chS.\Slabte} %xref
$\PD_X(\euler(\nu))=\sum^{m-1}_{n=0} (-1)^n \sum_{j \in J_n}
\{g_{j,n}\} \in H_1(X; \Z)$.
\endproclaim

\demo{Proof}
Using Definition C$_1$ of $\chi_1(X)(\gamma)$ given in
\cite{\HEC}, %bibref
and applying
Theorems 2.5 and 3.1 of \cite{\TTTF}, %bibref
we get
$$
\chi_1(X)(\gamma) = \left(
\sum^{m-1}_{n=0} (-1)^{n+1} \sum_{j \in J_n}
(|H_j|-1)\{g_{j,n}\} \right)
-
\PD_X(\euler(\nu)).
$$
Comparing this with the formula in the proof of
Theorem {\chS.\Slabtb}  %xref
yields the claimed result.\hfill{ }\qed
\enddemo

\remark{Remarks}
\roster
\item
Here is an informal explanation of the proof of
Theorem {\chS.\Slabte}. %xref
Assume that
\linebreak
$F^\gamma\:X \x I \ra X$ has no fixed points in
$X \x (0,\epsilon]$ where $\epsilon > 0$.
Moving the $[1-\epsilon,1]$ portion to $[-\epsilon,0]$,
let $\bar F^\gamma\: X \x [-\epsilon, 1-\epsilon] \ra X$ be
the reorganized homotopy.
The contribution of $\bar F^\gamma|_{X \x [-\epsilon, \epsilon]}$ to
$\chi_1(X)(\gamma)$ is $-\PD_X(\euler(\nu))$; this is 
Theorem 3.1 of \cite{\TTTF}. %bibref
The contribution of $\bar F^\gamma|_{X \x [\epsilon, 1-\epsilon]}$ to
$\chi_1(X)(\gamma)$ is the summation term; this is 
Theorem 2.5 of \cite{\TTTF}. %bibref

\item
In
Theorem {\chS.\Slabte} %xref
the right side is clearly independent of the choice of orientation
for $X$; however, so is the left side because
$\euler(-\nu) \cap [-X] = \euler(\nu) \cap [X]$.
\endroster
\endremark

Anticipating \S \chT, %xref
we apply Theorem {\chS.\Slabte} %xref
in the case $X$ is a compact connected
$3$--dimensional Seifert fibered space.
We will call a Seifert fibered space $X$ {\it admissible} if $X$
is oriented, the quotient surface $\Sigma$ is oriented,
and $X$ is not one of the special cases: $\Sigma = S^2$ with one
or two exceptional fibers,
or $\Sigma = D^2$ with one exceptional fiber.
Theorem {\chS.\Slabte}     %xref
and
Theorem {\chT.\Tlabpb}(1) %xref
yield the following formula:

\proclaim{Theorem \chS.\Slabtg} %xref
Let $X$ be an admissible $3$--dimensional Seifert fibered space
with quotient surface $\Sigma$ and $r$ exceptional fibers.
Impose the compatible orientation on the fibers
and on the normal bundle $\nu$ (to the Seifert fibering).
Let
$\{\gamma_0\} \in H_1(X;\Z)$ be the homology class of an ordinary
fiber, and let $\{g_j\}$ be the homology class of the
$j$--th singular fiber traversed once in the positive direction.
Then the Poincar\'e dual of the Euler class
of the oriented bundle $\nu$ is
$$
\PD(\euler(\nu)) = (\chi(\Sigma) - r)\{\gamma_0\} +
\sum_{j=1}^r \{g_j\}
~\in H_1(X; \Z).\hfill{ }\qed
$$
\endproclaim

%%%%%%%%%%%%%%%%%%%%%%%%%%
%%%   sec05.tex           
%%%%%%%%%%%%%%%%%%%%%%%%%%

\heading
\S \chV. $T^2$--actions
\endheading

In this section we show that our invariants vanish
for $S^1$--actions which extend to
\linebreak
$T^2$--actions
with finite isotropy.
The main results are
Theorem {\chV.\Vlabtb}, %xref
Theorem {\chV.\Vlabtc} %xref
and
Corollary {\chV.\Vlabca}. %xref

Let $T^2 \equiv S^1 \x S^1$ denote the $2$--torus viewed as
a compact Lie group.
A {\it $T^2$--space} is a topological space $X$ together
with a continuous left $T^2$--action $T^2 \x X \ra X$.
Given any circle subgroup $S^1 \subset T^2$, the restriction
of  $T^2 \x X \ra X$ to $S^1 \x X \ra X$ yields a circle
action.
Define $\Gamma^{T^2}_X \subset \Gamma_X$ to be the subgroup
generated by the fundamental elements
(see \S \chS) %xref
associated with circle actions obtained in this manner.
Note that $\Gamma^{T^2}_X$ is the image of the homomorphism
$\pi_1(T^2,1) \ra \Gamma_X$ induced by the map
$T^2 \ra X^X$ which is adjoint to the action map $T^2 \x X \ra X$.
%Hence $\Gamma^{T^2}_X$ is generated by fundamental elements
%$\gamma_1$ and $\gamma_2$ corresponding to the circle subgroups
%$S^1 \x \{1\}$ and $\{1\} \x S^1$ respectively.

Let $X$ be a $T^2$--CW complex
(the notion of a $K$--CW complex, for any compact Lie group $K$
is defined in
\cite{I$_1$}; %bibref
see Definitions {\chS.\Slabda} %xref
and {\chS.\Slabdb} %xref
for the case $K=S^1$).
For each
$T^2$--$n$--cell, $c^n_j \subset X$, 
there is a closed subgroup $H_j \subset T^2$
and a $T^2$--map
$$
f_j\:(T^2/H_j \x D^n, T^2/H_j \x S^{n-1}) \ra (c^n_j, \del c^n_j)
$$
such that $f_j(T^2/H_j \x D^n) = c^n_j$ and $f_j$ maps
$T^2/H_j \x D^n  ~-~ T^2/H_j \x S^{n-1}$ homeomorphically onto
$c^n_j - \del c^n_j$.
The subgroup $H_j$ is called the {\it isotropy group} of $c^n_j$
and the restriction of $f_j$ to $T^2/H_j \x S^{n-1}$ is called
the {\it attaching map} of $c^n_j$.
Let $\F = \{X_n\}$ be the filtration of $X$ by the
$T^2$--$n$--skeleta
($X_n$ is the union of the
$T^2$--$k$--cells of $X$, $k \leq n$).
We describe a cell structure on $\tilde T^2 \equiv \R^2$
which will be used to give
$(X,\F)$ the structure of a filtered cell complex.
Let $p\:\R^2 \ra T^2$ be the universal covering projection,
$p(t_1, t_2) = (e^{2\pi \imath t_1}, e^{2\pi \imath t_2})$.
The group of covering translations is $\Z^2$ acting by
$(m,n)\cdot(t_1,t_2) = (t_1+m, t_2+n)$ where
$(m,n) \in \Z^2$ and $(t_1,t_2) \in \R^2$.
Let $I = [0,1]$.
The cells
$\tE^2 = I \x I$,
$\tE^1_1 = I \x \{0\}$,
$\tE^1_2 = \{0\} \x I$,
$\tE^0 = \{(0,0)\}$
together with their translates give $\R^2$ a CW structure;
their images
$E^2   = p(\tE^2)$,
$E^1_1 = p(\tE^1_1)$,
$E^1_2 = p(\tE^1_2)$,
$E^0   = p(\tE^0)$
give $T^2$ a CW structure.
Let $(D^n,S^{n-1})$ have the standard relative CW structure
(with a single $n$--cell).
The space $T^2/H_j$ is given a CW structure which depends
on the dimension of the subgroup $H_j$:
\roster
\item
If $\dim H_j =2$, $H_j =T^2$ and so $T^2/H_j$ is a single point.

\item
If $\dim H_j=1$, choose an orientation preserving isomorphism 
of Lie groups 
\linebreak
$S^1\cong T^2/H_j$.
Give $S^1$ the CW structure consisting of one $0$--cell and
one \hbox{$1$--cell} and give $T^2/H_j$ the CW structure
induced by the chosen isomorphism.

\item
If $\dim H_j=0$, choose an orientation preserving isomorphism
of Lie groups
\linebreak
$T^2 \cong T^2/H_j$.
Give $T^2$ the CW structure described above
and give $T^2/H_j$ the CW structure induced by the chosen isomorphism.

\endroster
The pair $(T^2/H_j \x D^n, T^2/H_j \x S^{n-1})$ is given the
product relative CW structure. 
The maps $f_j$ determine a relative CW structure on
$(X_n, X_{n-1})$ thus realizing $(X,\F)$ as a filtered cell complex.

By Theorem 7.1 and Proposition 6.1 of \cite{I$_2$} %bibref
(where the case of a general compact Lie group action is treated),
we have:

\proclaim{Proposition {\chV.\Vlabta}} %xref
Any smooth $T^2$--manifold $M$ (with or without boundary)
has a
\linebreak
``$T^2$--triangulation'' and thus a
$T^2$--CW complex structure (which is finite if $M$ is compact)
such that the attaching maps of the $T^2$--cells are
piecewise linear.
In particular,
a $T^2$--triangulation gives rise to a polyhedral
filtered cell complex.\hfill{ }\qed
\endproclaim

Given a pair of integers $(a,b)$,
define $\tilde F^{(a,b)}\:\R^2 \x I \ra \R^2$ by
$$
\tilde F^{(a,b)}((t_1,t_2),t) =(t_1 + ta, t_2 + tb)
\qquad
(t_1,t_2)\in \R^2, ~t\in I.
$$
Note that $\tilde F^{(a,b)}$ descends to a map
$F^{(a,b)}\:T^2 \x I \ra T^2$ and
any map $F\:T^2 \x I \ra T^2$ with $F_0 = F_1 = \id$
is homotopic rel $T^2 \x \{0,1\}$ to $F^{(a,b)}$ for
some unique $(a,b)$.
This establishes an isomorphism of
groups $\Gamma_{T^2} @>{\cong}>> \Z \x \Z$ where
$\gamma \in \Gamma_{T^2}$, represented by
$F^{(a,b)}\:T^2 \x I \ra T^2$, is mapped to $(a,b)$.
The map $F^{(a,b)}$ is not necessarily cellular with respect
to the cell structure we have imposed on $T^2$;
however, $F^{(a,b)}$is homotopic  rel $T^2 \x \{0,1\}$ to the
cellular map $G^{(a,b)}\:T^2 \x I \ra T^2$ whose
lift $\tilde G^{(a,b)}\:\R^2 \x I \ra \R^2$ is given by
$$
\tilde G^{(a,b)}((t_1, t_2),t) =
\cases
(t_1, t_2 + 2tb)         &\text{if $0 \leq t \leq \tfrac{1}{2}$}  \\              
(t_1 + (2t-1)a, t_2 + b) &\text{if $\tfrac{1}{2} \leq t \leq 1$.} 
\endcases
$$
Let $H\equiv \pi_1(T^2,(1,1)) \cong \Z \x \Z$.
The {\it standard generators} of $H$ are the elements
$x_1$, $x_2$ represented respectively by
$I \ra S^1 \x S^1$, $t \mapsto (e^{2 \pi \imath t},1)$
and
$I \ra S^1 \x S^1$, $t \mapsto (1,e^{2 \pi \imath t})$.
For any $x \in H$ and integer $m$ define $x^{[m]} \in \Z H$ by
$$
x^{[m]} = 
\cases
(1 + x + \cdots + x^{m-1})            &\text{ if $m > 0$}\\
0                                     &\text{ if $m = 0$}\\
(-x^{-1} - x^{-2} - \cdots -x^{m})    &\text{ if $m<0$.}
\endcases
$$
Then it straightforward to show:

\proclaim{Lemma \chV.\Vlabma}  %xref
The chain homotopy
$\tilde D^{(a,b)}_*\:C_*(\tilde T^2) \ra C_{*+1}(\tilde T^2)$
associated to $G^{(a,b)}$ is given by
$$
\align
\tilde D^{(a,b)}_0(\tilde E^0) &= \tilde E^1_1 (x_1^{-1})^{[a]}x^{-b}_2 + 
                          \tilde E^1_2 (x_2^{-1})^{[b]}   \\
\tilde D^{(a,b)}_1(\tilde E^1_1) &= -\tilde E^2 (x_2^{-1})^{[b]}  \\
\tilde D^{(a,b)}_1(\tilde E^1_2) &= \tilde E^2 (x_1^{-1})^{[a]}x^{-b}_2.\hfill{ }\qed
\endalign
$$
\endproclaim

\proclaim{Theorem \chV.\Vlabtb}  %xref
Let $X$ be a finite connected $T^2$--CW complex such that all the
isotropy groups $H_j$ are finite.
Then $\tChi_1(X;\F)(\gamma)=0$ for all $\gamma \in \Gamma^{T^2}_X$.
\endproclaim

\demo{Proof}
For each $T^2$--$n$--cell $c^n_j$ with associated map
$f_j\:T^2/H_j \x D^n \ra c^n_j$ define
$$
g_{j,i,n} \equiv (f_j\circ (h,k))_{\#}(x_i) \in G \equiv \pi_1(X,v),
\qquad  i=1,2
$$
where $k\:T^2 \ra D^n$ is the constant map
(taking $T^2$ to the basepoint of $D^n$),
\linebreak
$h\:T^2 @>{\cong}>>T^2/H_j$ is the isomorphism
used to endow $T^2/H_j$ with its CW structure (note that
$H_j$ is finite by hypothesis) and $x_i$, $i=1,2$, are
the standard generators of $\pi_1(T^2,(1,1))$.
Let $e^{n+2}_j$, $e^{n+1}_{j,i}$, $i=1,2$, and $e^n_j$ be
the cells given by
$e^{n+2}_j = f_j \circ (h \x \id) (E^2 \x D^n)$,
$e^{n+1}_{j,i} = f_j \circ (h \x \id) (E^1_i \x D^n)$, $i=1,2$,
and $e^n_j = f_j \circ (h \x \id) (E^0 \x D^n)$ and
let
$\tilde e^{n+2}_j$, $\tilde e^{n+1}_{j,i}$, $i=1,2$,
and $\tilde e^n_j$ be corresponding lifts of these
cells to $\tilde X$.
The relative cellular chain complex $C_*(\tilde X_n,\tilde X_{n-1})$
is zero in degrees other than $n$, $n+1$, $n+2$:
$$
0 \ra C_{n+2}(\tilde X_n,\tilde X_{n-1})
@>{\tilde\partial^{n,n-1}_{n+2}}>>
C_{n+1}(\tilde X_n,\tilde X_{n-1})
@>{\tilde\partial^{n,n-1}_{n+1}}>>
C_n(\tilde X_n,\tilde X_{n-1}) \ra 0
$$
where the boundary operators are given by
$$
\align
\tilde\partial^{n,n-1}_{n+2}(\tilde e^{n+2}_j)      &=
(-1)^n \left(
\tilde e^{n+1}_{j,1}(1 - g^{-1}_{j,2,n}) -
\tilde e^{n+1}_{j,2}(1 - g^{-1}_{j,1,n})
\right) \\
\tilde\partial^{n,n-1}_{n+1}(\tilde e^{n+1}_{j,1})  &=
(-1)^{n+1} \tilde e^n_j(1 - g^{-1}_{j,1,n}) \\
\tilde\partial^{n,n-1}_{n+1}(\tilde e^{n+1}_{j,2})  &=
(-1)^{n+1} \tilde e^n_j(1 - g^{-1}_{j,2,n}). 
\endalign
$$

Let $\gamma \in \Gamma^{T^2}_X$ be the fundamental element
corresponding to a circle subgroup
\linebreak
$j\:S^1 \hookrightarrow T^2$,
i.e., if $\alpha\:T^2 \x X \ra X$ is the action map then
$\gamma$ is represented by $F^\gamma\:X \x I \ra X$,
$F^\gamma(x,t) =\alpha(j(e^{2 \pi \imath t}), x)$.
Clearly, $F^\gamma$ preserves the filtration $\F$;
however,
the maps $F^{\gamma,n,n-1}\:(X_n, X_{n-1}) \x I \ra (X_n, X_{n-1})$
are not necessarily cellular.
Using Cellular Approximation,
Lemma {\chV.\Vlabma} %xref
and the discussion preceding it,
we can replace $F^\gamma$ by a filtration preserving
map $G^\gamma$ such that the maps
$G^{\gamma,n,n-1}$ are cellular and the corresponding chain
homotopies, $\tilde D^{\gamma,n,n-1}_*$, are given by:
$$
\align
\tilde D^{\gamma,n,n-1}_n(\tilde e^n_j)
  &= \tilde e^{n+1}_{j,1} (g_{j,1,n}^{-1})^{[a_j]}g_{j,2,n}^{-b_j} + 
                          \tilde e^{n+1}_{j,2} (g_{j,2,n}^{-1})^{[b_j]}   \\
\tilde D^{\gamma,n,n-1}_{n+1}(\tilde e^{n+1}_{j,1})
  &= -\tilde e^{n+2}_j (g_{j,2,n}^{-1})^{[b_j]}  \\
\tilde D^{\gamma,n,n-1}_{n+1}(\tilde e^{n+1}_{j,2}) 
  &= \tilde e^{n+2}_j (g_{j,1,n}^{-1})^{[a_j]}g_{j,2,n}^{-b_j}
\endalign
$$
where $(a_j,b_j)$ is a pair of integers determined by the circle
subgroup $j\:S^1 \hookrightarrow T^2$ and the isotropy group $H_j$.
We write
$\tilde \partial^{n,n-1} = \bigoplus_k \tilde \partial^{n,n-1}_k$,
$\tilde D^{\gamma,n,n-1} = \bigoplus_k (-1)^{k+1} \tilde D^{\gamma,n,n-1}_k$
(viewed as matrices over $\Z G$;
see \S \chN.). %xref
Then
$$
\align
\trace(\tilde \partial^{n,n-1} \otimes \tilde D^{\gamma,n,n-1}) &= \sum_j \bigg(
   (-1)^n  (1 - g^{-1}_{j,2,n})  \otimes   -(-1)^{n+1} (g_{j,2,n}^{-1})^{[b_j]} \\
~&+~
-(-1)^n  (1 - g^{-1}_{j,1,n})  \otimes   (-1)^{n+1}(g_{j,1,n}^{-1})^{[a_j]}g_{j,2,n}^{-b_j} \\
~&+~
(-1)^{n+1} (1 - g^{-1}_{j,1,n})  \otimes   (-1)^n (g_{j,1,n}^{-1})^{[a_j]}g_{j,2,n}^{-b_j}\\
~&+~
(-1)^{n+1} (1 - g^{-1}_{j,2,n})  \otimes   (-1)^n  (g_{j,2,n}^{-1})^{[b_j]} \bigg).
\endalign
$$
By an obvious cancellation of terms with opposite sign on the right side of
the above expression, we conclude
$\trace(\tilde \partial^{n,n-1} \otimes \tilde D^{\gamma,n,n-1}) =0$.
Hence
$$
\tChi_1(X;\F)(\gamma) = \sum_n T_1(\tilde \partial^{n,n-1} \otimes \tilde D^{\gamma,n,n-1})
= \sum_n 0 = 0.\hfill{ }\qed
$$
\enddemo

Applying
Theorem \chF.\Flabta, %xref
we deduce:

\proclaim{Theorem \chV.\Vlabtc}  %xref
Suppose $X$ is a connected finite $T^2$--CW complex
such that all the isotropy groups are finite and
that $X$ is polyhedral as a filtered cell complex.
Let $i\:\Gamma^{T^2}_X \ra \Gamma_X$ be the inclusion.
Then
$\tChi''_1(X) \circ i = 0$
and
$i^*\tchi'_1(X) =0$.\hfill{ }\qed 
\endproclaim

By Theorem {\chS.\Slabta},  %xref
we obtain the following ``vanishing theorem'' for
the $S^1$--Euler characteristic of
{\S \chS} %xref
in the presence of a $T^2$--action with finite isotropy:

\proclaim{Corollary \chV.\Vlabca}  %xref
Suppose $X$ is a connected finite $S^1$--CW complex
which is polyhedral as a filtered cell complex
and there exists a $T^2$--action on $X$ with finite
isotropy groups such that $X$ has a $T^2$--CW complex
structure and the given $S^1$--action is obtained
by restriction of the $T^2$--action to a circle subgroup.
Then 
$\pi''(\tchi_{S^1}(X))=0$
and the derivation
$\langle \gamma \rangle \ra  HH_1(\Z G)'$ which takes  $\gamma$
to $\pi'(\tchi_{S^1}(X))$ is inner
where $\gamma \in \Gamma_X$ is the fundamental element
of the given $S^1$--action.\hfill{ }\qed
\endproclaim

In particular, the  $S^1$--Euler characteristic 
can be viewed as an obstruction
to the existence of a $T^2$--action with finite isotropy.

%%%%%%%%%%%%%%%%%%%%%%%%%%
%%%   sec06.tex           
%%%%%%%%%%%%%%%%%%%%%%%%%%

\heading
\S \chT.  Computations for $3$-Dimensional Seifert Fibered Spaces
\endheading

In this section we compute $\tchi_{S^1}(X)$ when $X$ is an
oriented $3$--dimensional Seifert fibered space
(Theorem \chT.\Tlabpb)
and discuss how much of that structure is detected by $\tchi_{S^1}(X)$
(Theorem \chT.\Tlabtc). %xref

Consider a smooth $S^1$--action without fixed points on a
compact connected oriented
$3$--manifold $X$.
The orbits of such an action
give $X$ the structure of a Seifert fibered space whose fibers
are consistently oriented,
and, conversely, any such Seifert fibered space comes from an
$S^1$--action
(see \cite{Sc, p.430}). %bibref
A Seifert fibering has only
finitely many singular fibers.
By Proposition \chS.\Slabtd, %xref
the $S^1$--space
$X$ can be given the structure of an
\linebreak
$S^1$--CW complex which is
polyhedral as a filtered cell complex so that the singular
fibers are $S^1$--$0$--cells.
Thus there are only finitely many
subgroups of $S^1$ which are orbit stabilizers.
Factoring out
the intersection of these stabilizers if necessary, we will
assume that the $S^1$--action is free away from the singular
fibers.
\footnote{
General references for the material on Seifert fibered
spaces used here are
\cite{B},        %bibref
\cite{J},        %bibref
\cite{JS} and    %bibref
\cite{Se}.}      %bibref

Our trace invariants are features of the $S^1$--action on $X$,
but a Seifert fibering, even when the fibers can be consistently
oriented, is a different type of structure on $X$.
In this section we explore
some of the ways in which our invariants give information
about the Seifert fibering.

The stabilizer $H_j$ is trivial unless $c^0_j$ is a singular fiber,
in which case we write
\linebreak
$\mu_j \equiv |H_j| > 1$.
Let $J'_0 \subset J_0$ be the indexing set for the singular
fibers
(as in \S \chS, %xref
$J_0$ indexes the $S^1$--$0$--cells).
We write $r\equiv |J'_0|$, the number of singular fibers, and,
abusing notation,
we sometimes identify $J'_0$ with $\{1,\ldots, r\}$.
We assume the basepoint $x_0$ lies in an ordinary fiber.
When making use of the homomorphism
$\eta_{\#}\:\Gamma_X \ra G \equiv \pi_1(X,x_0)$,
we write
$\gamma_0 = \eta_{\#}(\gamma)$, $\alpha_0 = \eta_{\#}(\alpha)$, etc.,
(recall that $\Gamma_X = \pi_1(X^X,\id_X)$ where $X^X$ is the function
space of self-maps of $X$ and that $\eta_{\#}$ is induced by evaluation
at the basepoint).
For $j \in J'_0$ we abbreviate $g_{j,0}$ to $g_j$.
As in
\S \chS, %xref
the $S^1$--action defines $F^\gamma\:X \x I \ra X$ where
$\gamma \in \Gamma_X$ is the fundamental element.
Note that $g_j^{\mu_j} = \gamma_0 \in Z(G)$ for each
$j \in J'_0$ (recall that $Z(G)$ denotes the center of $G$).

The quotient space $\Sigma \equiv X/S^1$ is a compact
connected surface.
A consistent orientation on the fibers of $X$ imposes an
orientation on $\Sigma$.
Conversely, if $\Sigma$ is
orientable then the fibers of $X$
can be consistently oriented.

We continue to assume $X$ and $\Sigma$ are oriented, and we recall
some facts about the fundamental group of $X$.

\subheading{Case 1. $\del X = \emptyset$}
The Seifert fibering is completely determined by
$\Sigma$, pairs of integers $(\mu_j, \beta_j)$ for $1\leq j \leq r$,
and an integer $b$.
Here,
the positive integers $\mu_j$ are as above.
For each singular fiber $c^0_j$,
there is a number $0 < \nu_j < \mu_j$
with $\text{gcd}(\nu_j,\mu_j)=1$ so that a fibered solid torus
neighborhood of $c^0_j$ is obtained from a standard fibered solid
torus by cutting and regluing using a $2 \pi\nu_j/ \mu_j$ twist;
and integers $\alpha_j$, $\beta_j$ are chosen so that
$\alpha_j \mu_j + \beta_j \nu_j =1$ and $0 < \beta_j < \mu_j$.
The integer $b$ measures the obstruction to constructing
a section $\Sigma \ra X_0$ where $X_0$ is obtained from $X$
by drilling out the singular fibers and filling the holes
with standard fibered solid tori.

Writing $\sigma$ for the genus of $\Sigma$,
a well known presentation of $G$
(see \cite{B}) %bibref
is:
$$
\langle
\gamma_0, a_1, \ldots, a_\sigma, b_1, \ldots, b_\sigma, 
c_1, \ldots, c_r ~|~
\text{$\gamma_0$ is central},~
c^{\mu_j}_j \gamma_0^{\beta_j}=1,~
\prod^\sigma_{i=1}[a_i,b_i] \prod^r_{j=1}c_j\gamma_0^b=1
\rangle
$$
This presentation is equivalent by Tietze transformations
$g_j \equiv c_j^{-\nu_j} \gamma_0^{\alpha_j}$
to the presentation:
$$
\langle
\gamma_0, a_1, \ldots, a_\sigma, b_1, \ldots, b_\sigma, 
g_1, \ldots, g_r ~|~
\text{$\gamma_0$ is central},~
g^{\mu_j}_j=\gamma_0,~
\prod^\sigma_{i=1}[a_i,b_i] \prod^r_{j=1}g^{-\beta_j}_j\gamma_0^b=1
\rangle
$$
One thinks of the $a_i$'s and $b_i$'s as generators of $\pi_1(\Sigma)$
and of $\gamma_0$ and
$g_j$
as having their previous meanings.

\subheading{Case 2. $\del X \neq \emptyset$}
Then $X$ is aspherical and $G$ is infinite.
We saw 
in {\S \chN} %xref
that in the aspherical case
$\Gamma_X \cong Z(G)$
so our
trace invariants only depend on the fundamental group.
\footnote{
While in the closed aspherical
manifold case the fundamental group determines the
Seifert fibered space up to homeomorphism,
this is not so when there is a non-empty boundary,
e.g\. if $\Sigma_1$ is a surface of genus $0$ with three
boundary components and
$\Sigma_2$ is a surface of genus $1$ with one
boundary component then the aspherical
Seifert fibered manifolds
$S^1 \x \Sigma_1$ and $S^1 \x \Sigma_2$ have the same fundamental
group but are not homeomorphic.}
In this case, 
after Tietze transformations as above,
$G$ has a presentation:
$$
\langle
\gamma_0, a_1, \ldots, a_\sigma, b_1, \ldots, b_\sigma, 
g_1, \ldots, g_r,  d_1, \ldots ,d_{m-1} ~|~
\text{$\gamma_0$ is central},~
g^{\mu_j}_j=\gamma_0
\rangle
$$
Here, $m > 0$ is the number
of boundary components of $\Sigma$,
while $\sigma$, $g_j$, $\gamma_0$ and $\mu_j$
are as before.
The numbers $\beta_j$ are absent from the
presentation because we ``solved'' for $d_m$ and removed it from
the set of generators as we removed the only relation
involving the $\beta_j$'s.

We will call a Seifert fibered space $X$ {\it admissible} if $X$
is oriented, $\Sigma$ is oriented,
and $X$ is not one of the special cases: $\Sigma = S^2$
and $r=1$ or $2$, or $\Sigma = D^2$ and $r=1$.

Recall that the Gottlieb subgroup of $G$, denoted by $\G (X)$,
is the image of the homomorphism $\eta_{\#}\:\Gamma_X \ra G$
induced by evaluation at the basepoint.
Our computations will depend on the following basic fact:

\proclaim{Proposition \chT.\Tlabpa} %xref
Let $X$ be admissible.
If $0 < i < \mu_j$ then $g_j^i \notin Z(G)$,
hence $g_j^i \notin \G (X)$.
If $0 < k < \mu_\ell$ and $(j,i) \neq (\ell, k)$ then
$g^i_j$ is not conjugate to $g^k_\ell$.
\endproclaim

\remark{Remark}
This is false when $X$ is not admissible.
\endremark

\demo{Proof of \chT.\Tlabpa} %xref
This is well known to experts so we confine ourselves to a
short sketch.

\flushpar
Case $1$: $\del X \neq \emptyset$; then the claimed facts can be
read off from the given presentation of $G$ by killing $\gamma_0$ and each
$a_i$, $b_i$ and $d_i$,
except for the non-admissible case $\Sigma=D^2$ and $r=1$.

\flushpar
Case $2$: $\del X = \emptyset$ and $r \geq 3$;
then one finds an appropriate triangle group as Fuchsian quotient;
the corresponding claims in the triangle group follow from
the geometry of its well known action on hyperbolic or
euclidean $2$--space or on the $2$--sphere
(see \cite{Si}). %bibref

\flushpar
Case $3$: $\del X = \emptyset$ and $r \leq 2$;
if $r = 0$ there is nothing to prove, so
assume $r = 1$ or $2$.
If $\Sigma$ has genus $\geq 2$, split $\Sigma$ along a
circle so that each part has genus at least $1$,
and arrange $X$ so that
the  singular fibers do not lie over this circle, and, if there are
two, that they lie one over each side.
This splits $G$ as a free product
with amalgamation, and a simple geometric argument establishes the claim.
A similar HNN argument handles the case where
$\Sigma$ has genus $1$.\hfill{ }\qed
\enddemo

We saw in {\S \chA} %xref
that $HH_*(\Z G) \cong \bigoplus_{C \in G_1} H_*(Z(g_C))$
by an isomorphism which is
canonical once a representative $g_C$ has been chosen
for each $C \in G_1$ (recall that $Z(g)$ denotes the centralizer
of $g \in G$ and that $G_1$ is the set of conjugacy classes of $G$).
In the present case
we choose  $\gamma_0^{-1}$ and $g^{-i}_j$ as representatives
of their conjugacy classes,
$j \in J'_0$, $0 < i < \mu_j$.
Since $\gamma_0 \in Z(G)$, we have
$G=Z(\gamma_0)$ and so
$H_1(G) \equiv H_1(Z(\gamma_0))  \cong HH_1(\Z G)_{C(\gamma_0)}$.
We remind the reader of the following notation convention:
given a subgroup $K \subset G$ and $g \in K$ we write
$\{g\} \in H_1(K) \equiv K_{\text{ab}}$ for the image of
$g$ under abelianization $K \ra K_{\text{ab}}$.

\proclaim{Theorem \chT.\Tlabpb} %xref
Let $X$ be admissible.
The components of $\tchi_{S^1}(X)$ in\newline
$\bigoplus_{C \in G_1} H_1(Z(g_C)) ~\cong~ HH_1(\Z G)$ are:
\roster
\item
The $C(\gamma_0)$--component is
$
(r - \chi(\Sigma))\{\gamma_0\} -
\sum_{j=1}^r \{g_j\} ~\in
H_1(G)
$.

\item
The $C(g^i_j)$--component is
$-\{g_j\} \in H_1(Z(g^i_j))$
for each $j \in J'_0$ and
$0 < i < \mu_j$.

\item
All other components are zero. \hfill{ }\qed
\endroster
\endproclaim 

\remark{Remark} In the notation of {\S \chN(A)}, %xref
(1) is the $HH_1(\Z G)'$ part and (2) is the $HH_1(\Z G)''$ part.
\endremark

\demo{Proof}
In canonical form,
$\tchi_{S^1}(X)$ is represented by the cycle
$$
\zeta ~=~ \sum_{n \geq 0} (-1)^{n+1} \sum_{j \in J_n - J'_0}
         \gamma_0 \otimes \gamma_0^{-1} \gamma_0^{-1} 
       ~-~ \sum_{j=1}^r \sum_{i=1}^{\mu_j}
         g_j \otimes g_j^{-1} g_j^{-i} 
       ~=~ \zeta' + \zeta''
$$
where
$$
\zeta' ~=~ \sum_{n \geq 0} (-1)^{n+1} \sum_{j \in J_n - J'_0}
         \gamma_0 \otimes \gamma_0^{-1} \gamma_0^{-1} 
       ~-~ \sum_{j \in J'_0} g_j \otimes g_j^{-1} \gamma_0^{-1}
$$
and
$$
\zeta'' ~=~ -\sum_{j=1}^r \sum_{i=1}^{\mu_j-1}
              g_j \otimes g_j^{-1}  g_j^{-i}. 
$$
By Proposition {\chT.\Tlabpa},  %xref
this decomposes $\zeta$ into a single central component and
various components whose markers are not in $\G (X)$.

For any $a \in G$, $a \otimes a^{-2}$ is homologous  to
$-a^{-1} \otimes 1 \equiv -a^{-1} \otimes a a^{-1}$.
Thus, writing ``$u \sim v$'' when $u$ and $v$ are homologous
chains in the Hochschild complex
(see \S \chA), %xref
we have
$$
\zeta' \sim (\chi(\Sigma)-r)  \gamma_0^{-1} \otimes \gamma_0 \gamma_0^{-1} 
          + \sum_{j=1}^r g_j^{-1} \otimes g_j \gamma_0^{-1}
$$
and
$$
\zeta'' \sim \theta'' \equiv
\sum_{j=1}^r \sum_{i=1}^{\mu_j-1} g_j^{-1} \otimes g_j  g_j^{-i}. 
$$

For $c \in Z(G)$ and $a, b \in G$, we have
(see \cite{\AJM, Lemma 6.13})  %bibref
$$
ab^{-1} \otimes ba^{-1}c ~\sim~ a \otimes a^{-1}c - b \otimes b^{-1}c.
$$
So $\zeta' \sim \theta' \equiv h \otimes h^{-1} \gamma_0^{-1}$
where $h = \gamma_0^{r - \chi(\Sigma)}g_1^{-1}\cdots g_r^{-1}$.
Thus $\theta \equiv \theta' + \theta''$ is a representative Hochschild
cycle for $\tchi_{S^1}(X)$ in which every term lies in a different
component $C_*(\Z G)_C$.
Collecting terms,
$$
\theta \equiv \theta' + \theta'' = h \otimes h^{-1} \gamma_0^{-1}
~+~ 
\sum_{j=1}^r \sum_{i=1}^{\mu_j-1} g_j^{-1} \otimes g_j  g_j^{-i}.
$$
where $h = \gamma_0^{r - \chi(\Sigma)}g_1^{-1}\cdots g_r^{-1}$.
The conclusion of the theorem follows from
Proposition {\chA.\Alabpa} %xref
applied to this expression for $\theta$.\hfill{ }\qed
\enddemo

By Theorem \chS.\Slabta(2), %xref
Theorem \chT.\Tlabpb(1) %xref
and
Theorem {\chG.\Glabtb} %xref
we have:

\proclaim{Proposition \chT.\Tlabtb} %xref
Let $X$ be admissible and suppose
$(\chi(\Sigma) -r)\{\gamma_0\} +
\sum^r_{j=1} \{g_j\} \neq 0 \in H_1(X; \Z)$.
Then the Gottlieb subgroup $\G(X)$ is cyclic and is generated
by $\gamma_0$.\hfill{ }\qed
\endproclaim

\proclaim{Corollary \chT.\Tlabca} %xref
If in addition $X$ is aspherical
then $Z(G)$ is infinite cyclic and is generated
by $\gamma_0$.
\endproclaim

\demo{Proof}
By a well known argument
(see \cite{CR, p.43}) %bibref
$\gamma_0$ is non-trivial when $X$ is aspherical.\hfill{ }\qed
\enddemo

When $\{\gamma_0\}$ is non-zero in $H_1(G)$
each $\{g_j\}$ is non-zero in
$H_1(Z(g^i_j))$ for $0 < i < \mu_j$ since $g^{\mu_j}_j = \gamma_0$.
Thus,
combining
Theorems {\chS.\Slabta} %xref
and {\chS.\Slabtb} %xref
and
Theorem {\chT.\Tlabpb}, %xref
we see what features of the Seifert fibering
are detected by the {\it topological} invariants
$\tChi''_1(X)$, $\tchi'_1(X)$ and
$\chi_1(X)$ of {\S \chN}: %xref

\proclaim{Proposition \chT.\Tlabpd}  %xref
If $\{\gamma_0\} \neq 0 \in H_1(G)$ 
and $X$ is admissible then:
\roster
\item
$\tChi''_1(X)(\gamma)$ detects the number, $r$, of singular
fibers,
the numbers $\mu_j$ and the conjugacy
classes of the singular fibers.

\item
$\chi'_1(X)(\gamma)
=(r-\chi(\Sigma))\{\gamma_0\} -
\sum^r_{j=1} \{g_j\} \in H_1(G)$.

\item
$\chi_1(X)(\gamma) = -\chi(\Sigma)\{\gamma_0\} \in H_1(G)$.\hfill{ }\qed
\endroster
\endproclaim

Next, we mention the connection with two classical invariants of Seifert
fibered spaces,
the Euler number and the orbifold Euler characteristic.
When $\del X = \emptyset$ the {\it Euler number} of the Seifert
fibered space $X$
(see \cite{Sc, \S 5}) %bibref
is defined to be the rational number\linebreak
$b - \sum^r_{j=1} \beta_j/\mu_j$.
By abelianizing the given presentations for $G$
we conclude a well known fact:

\proclaim{Lemma \chT.\Tlabpc}  %xref
$\{\gamma_0\} \in H_1(G)$
has infinite order if and only if
either $\del X \neq \emptyset$ or 
$\del X =\emptyset$ and the Euler number of $X$ is
zero.\hfill{ }\qed
\endproclaim

Combining
Lemma {\chT.\Tlabpc} %xref
with
Theorems {\chT.\Tlabpb}
and {\chS.\Slabtb} %xref
we obtain:

\proclaim{Theorem \chT.\Tlabtc}  %xref
Suppose that $X$ is admissible
and either $\del X \neq \emptyset$ or the Euler number is zero.
Then $\tchi_{S^1}(X)$
determines the following
features of the Seifert fibering:
$\chi(\Sigma)$,
the number of singular fibers $r$,
the integers $\mu_1, \ldots, \mu_r$,
and the conjugacy classes in $G$
represented by the singular fibers.\hfill{ }\qed
\endproclaim

The {\it orbifold Euler characteristic} of the quotient surface $\Sigma$
(viewed as an orbifold) is the rational number
$\chi_V(\Sigma) \equiv \chi(\Sigma) + \sum_{j=1}^r (\tfrac{1}{\mu_j}-1)$.
Denote the image of a homology class $z \in H_1(X; \Z)$ in $H_1(X; \Q)$
by $z_{\Q}$.

\proclaim{Proposition \chT.\Tlabpe} %xref
$\pi'(\tchi_{S^1}(X))_{\Q} = -\chi_V(X)\{\gamma_0\}_{\Q}$.
\endproclaim

\demo{Proof}
We have $\{g_j\}_{\Q}= \tfrac{1}{\mu_j}\{\gamma_0\}_{\Q}$,
and so, 
by Theorem {\chT.\Tlabpb}(1) %xref
and Theorem {\chS.\Slabtg}: %xref
$$
\align
-\pi'(\tchi_{S^1}(X))_{\Q} =
\PD_X(\euler(\nu))_{\Q} &=
(\chi(\Sigma)-r) \{\gamma_0\}_{\Q} + 
\sum_{j=1}^r \tfrac{1}{\mu_j}\{\gamma_0\}_{\Q} \\
&= \bigg(\chi(\Sigma) +
\sum_{j=1}^r (\tfrac{1}{\mu_j}-1)\bigg)\{\gamma_0\}_{\Q}. \hfill{ }\qed
\endalign
$$
\enddemo
\remark{Remark}
In particular, $\tchi_{S^1}(X)$ determines $\chi_V(\Sigma)$ whenever
$\{\gamma_0\}$ has infinite order in $H_1(G)$.
\endremark

   From Proposition {\chT.\Tlabpc}  %xref
we conclude an important theoretical fact
about the connection between
the invariants studied in this paper and the group $K_1(\Z G)$.
Whenever one considers an invariant lying in $HH_1(\Z G)$
it is natural to ask if it
is the image of a richer invariant lying in $K_1(\Z G)$ under the Dennis
Trace map $\DT\:K_1(\Z G) \ra HH_1(\Z G)$; we discuss this in a number of
related contexts in \cite{GN$_4$}. %xref
Proposition {\chT.\Tlabpc}  %xref
shows that in the
case of $\tchi_1(X)$,
with $X$ a Seifert fiber space, the answer is often negative.
In more detail: 
$\tChi''_1(X)(\gamma)$ is represented by
$\sum_{j=1}^r \sum_{i=1}^{\mu_j-1} g_j^{-1} \otimes g_j  g_j^{-i}$.
If the Euler number
is zero and $X$ is aspherical and
without boundary,  then,
by Proposition \chT.\Tlabpc, %xref
$\{g_j\} \neq 0 \in H_1(Z(g_j))$
because $\mu_j \{g_j\} = {\gamma_0} \neq 0 \in H_1(G)$.
Thus if $r > 0$ there is a non-zero component of
$\tChi_1(X)(\gamma)$ marked by a non-central element.
Yet $\Wh_1(G) = 0$
(\cite{W})  %bibref
and so it follows from
Proposition {\chA.\Alabpc} %xref
that $\tChi_1(X)(\gamma)$ does not lie in the image of the Dennis trace
$\DT\: K_1(\Z G) \ra HH_1(\Z G)$.
Indeed $\DT$ extends to a $\Gamma_X$--module homomorphism
$\DT'\: \Z \Gamma_X \otimes K_1(\Z G) \ra HH_1(\Z G)$
and $\tchi_1(X)$ is not in the image of
$\DT'_*\:H^1(\Gamma_X, \Z \Gamma_X \otimes K_1(\Z G))
\ra H^1(\Gamma_X, HH_1(\Z G))$.

%%%%%%%%%%%%%%%%%%%%%%%%%%
%%%   secapA.tex          
%%%%%%%%%%%%%%%%%%%%%%%%%%

\heading
Appendix \chG.  Some consequences of $\tchi'_1(X) \neq 0$
\endheading

In {\S \chN} %xref
we defined a cohomomolgy class
$\tchi'_1(X) \in H^1(\Gamma_X, H_1(G) \otimes {\Cal G}(X))$
which is a homotopy invariant of the space $X$
(see Definition \chN.\Nlabdb). %xref
Recall that $\Gamma_X = \pi_1(X^X,\id_X)$ where $X^X$ is the function
space of self-maps of $X$ and ${\Cal G}(X)$ is the Gottlieb subgroup
of $G \equiv \pi_1(X,x_0)$;
see \S \chN. %xref
In this Appendix we examine some group theoretic and topological
consequences of the hypothesis $\tchi'_1(X) \neq 0$.

\proclaim{Theorem \chG.\Glabta}  %xref
Suppose $X$ is a finite aspherical CW complex
and  $\tchi'_1(X) \neq 0$.
Then $\Gamma_X \cong Z(G)$ is infinite cyclic.
\endproclaim

\demo{Proof}
Since $X$ is aspherical,
by \cite{Got} %bibref
the evaluation at the basepoint $\eta\:X^X \ra X$ induces an
isomorphism $\Gamma_X \cong Z(G)$.
By hypothesis $\tchi'_1(X) \neq 0$;
in particular, the cohomology group
$H^1(Z(G), H_1(G) \otimes_{\Z} \Z Z(G))$
is non-trivial.
Since $G$ is torsion free,
Theorem 6.10 of \cite{DD, Ch.IV}, %bibref
implies that $Z(G)$ is either infinite cyclic or
a non-trivial free product.
The group $Z(G)$ is abelian and so cannot be
a non-trivial free product.\hfill{ }\qed
\enddemo

In practice, we often have more information about $\tchi_1'(X)$
(e.g\. an explicit calculation of $\tChi'_1(X)(\gamma)$ for
certain $\gamma \in \Gamma_X$),
leading to a sharper conclusion:

\proclaim{Theorem \chG.\Glabtb}  %xref
Suppose $X$ is a finite CW complex,
$\gamma \in \Gamma_X$,
and  $\tchi'_1(X)$ is represented by a derivation
$\Delta\:\Gamma_X \ra HH_1(\Z G)'$ such that $\Delta(\gamma)$ has
a non-zero entry in exactly one component.
Then the Gottlieb subgroup, $\G(X)$, is cyclic and is
generated by $\eta_{\#}(\gamma)$.
\endproclaim

\remark{Remark}
Let $\mu$ be the composite homomorphism
$$
HH_1(\Z G)' \cong H_1(G) \otimes \Z \G(X)  @>{\id \otimes \vep}>> H_1(G)
$$
where $\vep\:\Z \G(X) \ra \Z$ is the augmentation.
The derivation $\Delta$ of
Theorem {\chG.\Glabtb}  %xref
is necessarily outer because $\mu(\Delta(\gamma)) \neq 0$ whereas
$\mu(J(\gamma))=0$ for any inner derivation $J$.
In particular  $\tchi'_1(X) \neq 0$;
compare
Theorem {\chG.\Glabta}.  %xref
\endremark

\demo{Proof of Theorem {\chG.\Glabtb}}  %xref
We may assume that $\Delta(\gamma)$ is concentrated in the
$\gamma^0$--component (i.e., the component corresponding
to the identity element).
Since $\Delta(\gamma) \neq 0$, $\gamma \neq 1$.
Let $\tau \in \Gamma_X$.
Assuming $\eta_{\#}(\tau) \neq 1$,
we will show $\eta_{\#}(\tau)$ is a power of $\eta_{\#}(\gamma)$.
since $\tau \gamma = \gamma \tau$,
we have $(1-\tau)\Delta(\gamma)=(1-\gamma)\Delta(\tau)$.
For any $\lambda \in HH_1(\Z G)'$,
the sum in $H_1(G)$ of the entries of $(1-\gamma)\lambda$
in a $\langle \gamma \rangle$--orbit of components of
$HH_1(\Z G)'$ is obviously zero.
Applying this observation with $\lambda = \Delta(\gamma)$,
we find that
$C(\eta_{\#}(\tau)^{-1}) = C(\eta_{\#}(\gamma)^{-k})$
for some $k$.
Hence $\eta_{\#}(\tau) = \eta_{\#}(\gamma)^k$ because
$\G(X) \equiv \eta_{\#}(\Gamma_X) \leq Z(G)$.\hfill{ }\qed
\enddemo

\proclaim{Corollary \chG.\Glabca}  %xref
If $X$ in
Theorem {\chG.\Glabtb}  %xref
is aspherical,
then $Z(G)$ is infinite cyclic and is generated by
$\gamma \in \Gamma_X \cong Z(G)$.\hfill{ }\qed
\endproclaim

The hypothesis in
Corollary {\chG.\Glabca}  %xref
that $\Delta(\gamma)$ be concentrated in one component is
necessary in the following sense:

\proclaim{Proposition \chG.\Glabpa}  %xref
Suppose $Z(G)$ is infinite cyclic with generator $\gamma$.
Let
$$
\Delta\:Z(G) @>{}>> HH_1(\Z G)' \cong 
 \bigoplus_{\gamma^n \in Z(G)} H_1(G)
$$
be a derivation.
Then $\Delta$ differs by an inner derivation
from a derivation $\Delta'$ such that
$\Delta'(\gamma)$ has a non-zero entry in at most one component.
\endproclaim

\demo{Proof}
A derivation 
$Z(G)\equiv\langle \gamma \rangle \ra HH_1(\Z G)'$
is freely determined by its value on $\gamma$.
Write $\Delta(\gamma) = \sum_i \Delta_i(\gamma)$ where $\Delta_i(\gamma)$
has the same $\gamma^i$--component, say $a_i$, as
$\Delta(\gamma)$ and all other components are zero.
Define a derivation $\Delta'_i$ so that 
the $\gamma^0$--component of $\Delta'_i(\gamma)$ is $a_i$
and all other components are zero.
Then
$
\Delta_i(\gamma) - \Delta'_i(\gamma) = (1 - \gamma^i) \Delta_i(\gamma)
= (1 - \gamma) u_i
$
where
$$
u_i = 
\cases
(1 + \gamma + \cdots + \gamma^{i-1})\Delta_i(\gamma) &\text{ if $i > 0$}\\
0  &\text{ if $i = 0$}\\
(-\gamma^{-1} - \gamma^{-2} - \cdots - \gamma^{i})\Delta_i(\gamma) &\text{ if $i<0$.}
\endcases
$$
Define $\Delta'(\gamma) = \sum_i \Delta'_i(\gamma)$.
Then $\Delta(\gamma) - \Delta'(\gamma) = (1 - \gamma)(\sum_i u_i)$ and
$\Delta'(\gamma)$ is concentrated in the $\gamma^0$--component.\hfill{ }\qed
\enddemo

\remark{Remark \chG.\Glabra} %xref
Recall that
$\chi_1(X) = \chi'_1(X) + \chi''_1(X)$
(see Definition \chN.\Nlabde). %xref
With a different and more
difficult proof,
Theorem 5.4 of \cite{\HEC} %bibref
provides the same conclusion as
Theorem {\chG.\Glabta}  %xref
under the hypotheses that $\chi_1(X; \Q) \neq 0$ and
$G$ has the ``Weak Bass Property over $\Q$''
(see \cite{\HEC, \S 5}). %bibref
Here, $\chi_1(X; \Q)\:\Gamma_X \ra H_1(X;\Q)$ is
the composite of $\chi_1(X)$ with the coefficient
homomorphism $H_1(X;\Z) \ra H_1(X;\Q)$.
In
Theorem {\chG.\Glabta},  %xref
the (possibly vacuous) $K$--theoretic hypothesis
about the Weak Bass Property is absent.
The hypothesis $\chi_1(X; \Q) \neq 0$
can be checked in many cases in view of the wide variety of equivalent
definitions of $\chi_1(X)$ given in \cite{\HEC}. %bibref
However, it can happen
that $\chi_1(X) = 0$ while $\tchi_1'(X) \neq 0$;
see Example 3.8 of  \cite{\HEC}. %bibref
\endremark

%%%%%%%%%%%%%%%%%%%%%%%%%%%%%%%%%%%%%%%%%%%%%%%%%%%%%%%%%%
%%   references list for:
%%   A Hochschild Homology Euler Characteristic for Circle Actions   
%%%%%%%%%%%%%%%%%%%%%%%%%%%%%%%%%%%%%%%%%%%%%%%%%%%%%%%%%%
%%   for use with AMSTeX version 2.1
%%%%%%%%%%%%%%%%%%%%%%%%%%%%%%%%%%%%%%%%%%%%%%%%%%%%%%%%%%

%%%%%%%%%%%%%%%%%%%%%%%%%%
%%%    refs.tex           
%%%%%%%%%%%%%%%%%%%%%%%%%%

\enddocument